\documentclass[11pt]{article}
\usepackage{color}
\usepackage{amsfonts}
\usepackage{amsmath}
\usepackage{amssymb}
\usepackage{latexsym}
\begin{document}
\newtheorem{theor}{Theorem}
\newtheorem{prop}{Proposition}
\newtheorem{lemma}{Lemma}
\newtheorem{coro}{Corollary}
\newtheorem{prof}{Proof}
\newtheorem{defi}{Definition}
\newtheorem{rem}{Remark}

\definecolor{dark}{gray}{.5}
\definecolor{light}{gray}{0.75}
\definecolor{midmid}{gray}{0.625}

\font\tenbb=msbm10
\font\sevenbb=msbm7
\font\fivebb=msbm5
\newfam\bbfam
\textfont\bbfam=\tenbb
\scriptfont\bbfam=\sevenbb
\scriptscriptfont\bbfam=\fivebb
\def\Bbb{\fam\bbfam\tenbb}


\newcommand{\begeq}[1]{\begin{equation} \label{#1}}
\newcommand{\fineq}{\end{equation}}

\newcommand{\ra}{\rightarrow}
\newcommand{\lra}{\longrightarrow}

\newcommand{\m}{\displaystyle}

\newcommand{\imp}{\Rightarrow}
\newcommand{\tend}{\rightarrow}
\newcommand{\g}{\gamma}
\newcommand{\G}{\Gamma}
\newcommand{\De}{\Delta}
\newcommand{\s}{\sigma}
\newcommand{\Si}{\Sigma}
\newcommand{\al}{\alpha}
\newcommand{\e}{\epsilon}
\newcommand{\de}{\delta}
\newcommand{\ka}{\kappa}
\newcommand{\la}{\lambda}
\newcommand{\ti}{\tilde}
\newcommand{\ov}{\overline}
\newcommand{\F}{{\cal F}}
\newcommand{\be}{\beta}
\newcommand{\ch}{\hat}
\newcommand{\un}{\mbox{\sl\em 1\hspace{-3pt}l}}


\def\bR{{\Bbb R}}                  
\def\bN{{\Bbb N}}                  
\def\bZ{{\Bbb Z}}                  
\def\bE{{\Bbb E}}                  
\def\bV{{\Bbb V}}                  
\def\bP{{\Bbb P}}                  
\def\bI{{\Bbb I}}                  

\def\tinf{{\rightarrow\infty}}     

\def\tts{{\textstyle}}

\def\Kd{{K^{(d)}}}
\def\fd{{f^{(d)}}}
\def\fnd{{f^{(d)}_n}}
\def\vn{{\sqrt{\frac{nh_n^{2d+1}}{a_n}}}}
\def\Lnxu{{\Lambda_{n,x}(u)}}
\def\Lxu{{\Lambda_x(u)}}


\begin{titlepage}
\title{The stochastic approximation method for the estimation of 
a multivariate probability density}
\author{Abdelkader Mokkadem \and Mariane Pelletier \and Yousri Slaoui}
\date{\small{
Universit\'e de Versailles-Saint-Quentin\\
D\'epartement de Math\'ematiques\\
45, Avenue des Etats-Unis\\
78035 Versailles Cedex\\
France\\}
(mokkadem,pelletier,slaoui)@math.uvsq.fr\\}
\maketitle

\begin{abstract}
We apply the stochastic approximation method to construct a large class 
of recursive kernel estimators of a probability density, including the one 
introduced by Hall and Patil (1994). We study the properties of 
these estimators and compare them with Rosenblatt's nonrecursive estimator. 
It turns out that, for pointwise estimation, 
it is preferable to use the nonrecursive Rosenblatt's kernel estimator 
rather than any recursive estimator. A contrario, for estimation by 
confidence intervals, it is better to use a recursive estimator rather than 
Rosenblatt's estimator. 
\end{abstract}
\vspace{0.5cm}

\noindent
{\bf Key words and phrases~:}\\
{Density estimation; Stochastic approximation algorithm.}\\
\vspace{0.5cm}

\end{titlepage}

\section{Introduction}

The advantage of recursive estimators on their nonrecursive version is that their 
update, from a sample of size $n$ to one of size $n+1$, requires considerably 
less computations. This property is particularly important in the framework of 
density estimation, since the number of points at which the function is estimated 
is usually very large. The first recursive version of Rosenblatt's 
kernel density estimator - and the most famous one - was introduced by 
Wolwerton and Wagner (1969), and was
widely studied; see among many others Yamato (1971), Davies (1973), Devroye
(1979), Wegman and Davies (1979) and Roussas (1992). 
Competing recursive estimators, which may be regarded as weighted versions of Wolwerton and Wagner's
estimator, were introduced and studied by Deheuvels (1973),  Wegman and Davies (1979) and Duflo (1997).
Recently, Hall and Patil (1994) defined a large class of weighted recursive estimators, 
including all the previous recursive estimators.
In this paper, we apply the stochastic approximation method to define a class of recursive kernel 
density estimators, which includes the one introduced by Hall and Patil (1994).\\

The most famous use of stochastic approximation algorithms in the framework of 
nonparametric statistics is the work of Kiefer and Wolfowitz (1952), who 
build up an algorithm which allows the approximation of the maximizer of a 
regression function. Their well-known algorithm was 
widely discussed and extended in many directions (see, among many others,
Blum (1954), 
Fabian (1967),
Kushner and Clark (1978),
Hall and Heyde (1980), 
Ruppert (1982),
Chen (1988), 
Spall (1988), 
Polyak and Tsybakov (1990),
Dippon and Renz (1997),
Spall (1997),  
Chen, Duncan and Pasik-Duncan (1999), 
Dippon (2003), and Mokkadem and Pelletier (2007a)). 
Stochastic approximation algorithms were also introduced by 
R\'ev\'esz (1973, 1977) to estimate a regression function, and by 
Tsybakov (1990) to approximate the mode of a probability density.
\\

Let us recall Robbins-Monro's scheme to construct approximation algorithms of search of 
the zero $z^*$ of an unknown function $h~:\bR\tend\bR$. First, 
$Z_0\in\bR$ is arbitrarily chosen, and then the sequence $(Z_n)$ is recursively 
defined by setting  
$$Z_n=Z_{n-1}+\g_nW_n$$
where $W_n$ is an ``observation'' of the function $h$ at the point $Z_{n-1}$, and where 
the stepsize $(\g_n)$ is a sequence of positive real numbers that goes to zero.\\

Let $X_1,\ldots ,X_n$ be independent, identically distributed 
$\bR^d$-valued random vectors, 
and let $f$ denote the probability density of $X_1$. To construct a stochastic 
algorithm, which approximates the function $f$ at a given point $x$, we define an 
algorithm of search of the zero of the function $h~: y\mapsto f(x)-y$. We thus 
proceed in the following way: (i) we set $f_0(x)\in\bR$; 
(ii) for all $n\geq 1$, we set  
\begin{eqnarray*}
f_n(x)=f_{n-1}(x)+\g_nW_n(x)
\end{eqnarray*}
where $W_n(x)$ is an ``observation'' of the function $h$ at the point 
$f_{n-1}(x)$. To define $W_n(x)$, we follow the approach of R\'ev\'esz (1973, 1977) and 
of Tsybakov (1990), and introduce a kernel $K$ (that is, a function satisfying  
$\int_{\bR^d} K(x)dx=1$) and a bandwidth $(h_n)$ (that is, a sequence of positive real 
numbers that goes to zero), and set
$W_n(x)=h_n^{-d}K(h_n^{-1}[x-X_n])-f_{n-1}(x)$. 
The stochastic approximation algorithm we introduce to recursively estimate the 
density $f$ at the point $x$ can thus be written as
\begin{equation}\label{eq:1}
f_n(x)=(1-\g_n)f_{n-1}(x)+\g_nh_n^{-d}K\left(\frac{x-X_n}{h_n}\right).
\end{equation}
Let $(w_n)$ be a positive sequence such that $\sum w_n=\infty$. When the 
stepsize $\left(\gamma_n\right)$ is chosen 
equal to $(w_n[\sum_{k=1}^nw_k]^{-1})$, the estimator $f_n$ defined 
by~\eqref{eq:1} can be rewritten as
\begin{equation}
\label{peter}
f_n(x)=\frac{1}{\sum_{k=1}^nw_k}\sum_{k=1}^nw_k\frac{1}{h_k^d}
K\left(\frac{x-X_k}{h_k}\right).
\end{equation}
The class of estimators defined by the stochastic approximation algorithm~\eqref{eq:1} 
thus includes the general class 
of recursive estimators expressed as~\eqref{peter}, and introduced in Hall and Patil 
(1994). In particular, the choice $(w_n)=1$ produces the estimator proposed by 
Wolverton and Wagner (1969), the choice $(w_n)=(h_n^{d/2})$ yields the estimator 
considered by Wegman and Davies (1979), and the choice 
$(w_n)=(h_n^{d})$ gives the estimator considered by 
Deheuvels (1973) and Duflo (1997).\\

The aim of this paper is the study of the properties of the recursive 
estimator defined by the stochastic approximation algorithm~\eqref{eq:1}, 
and its comparison with the wellknown nonrecursive kernel density estimator 
introduced by Rosenblatt (1956) (\emph{see} also Parzen (1962)), and defined as
\begin{eqnarray}\label{eq:4}
\tilde{f}_n\left(x\right)=\frac{1}{nh_n^d}\sum_{k=1}^nK\left(\frac{x-X_k}{h_n}\right).
\end{eqnarray}
\\


We first compute the bias and the variance of the recursive 
estimator $f_n$ defined by~\eqref{eq:1}. It turns out that they heavily depend on 
the choice of the stepsize $\left(\gamma_n\right)$. In particular, for a given bandwidth, 
there is a trade-off in the choice of $\left(\gamma_n\right)$ 
between minimizing either the bias or the variance of $f_n$. To determine the 
optimal choice of stepsize, we consider two points of view: 
pointwise estimation and estimation by confidence intervals.\\

From the pointwise estimation point of view, the criteria we consider to find the 
optimal stepsize is minimizing the mean squared error (MSE) or the integrated mean 
squared error (MISE). 
We display a set of stepsizes $\left(\gamma_n\right)$ minimizing 
the MSE or the MISE of the estimator $f_n$ defined by~\eqref{eq:1};
we show in particular that the sequence $\left(\gamma_n\right)=\left(n^{-1}\right)$ 
belongs to this set.  
The recursive estimator introduced by Wolverton and Wagner (1969) thus 
belongs to the subclass of recursive kernel estimators which have a minimum MSE or 
MISE (thanks to an adequate choice of the bandwidth, \emph{see} 
Section~\ref{section 2.1}). Let us underline that these 
minimum MSE and MISE 
are larger than those obtained for Rosenblatt's nonrecursive estimator $\tilde f_n$.
Thus, for pointwise estimation and when rapid updating is not such important, it 
is preferable to use  
Rosenblatt's estimator rather than any recursive estimator defined by the 
stochastic approximation algorithm \eqref{eq:1}. Let us also mention that Hall and 
Patil (1994) introduce a class of on-line estimators, constructed from the class of 
the recursive estimators defined in \eqref{peter}; their on-line estimators are not 
recursive any more, but updating them requires much less operations than updating 
Rosenblatt's estimator, and their MSE and MISE are smaller than those of the 
recursive estimators \eqref{peter}.\\

Let us now consider the estimation from confidence interval point of view. Hall (1992) shows that, to minimize the coverage error of probability density confidence intervals, avoiding bias estimation by a slight undersmoothing is more efficient than explicit bias correction. In the framework of undersmoothing, minimizing the MSE comes down to minimizing the variance. We thus display a set of stepsizes $\left(\gamma_n\right)$ minimizing 
the variance of  $f_n$; we show in particular that, when the bandwidth $(h_n)$ 
varies regularly with exponent $-a$, the sequence 
$\left(\gamma_n\right)=\left([1-ad]n^{-1}\right)$ belongs to this set. 
Let us underline that the variance of the estimator $f_n$ defined with 
this stepsize is smaller than that of Rosenblatt's estimator. 
Consequently, even in the case when the on-line aspect is not quite important, it is preferable 
to use recursive estimators to construct confidence intervals. The simulation results 
given in Section \ref{section simul} are corroborating these theoritical results. 
\\


To complete the study of the asymptotic properties of the recursive estimator $f_n$, 
we give its pointwise strong convergence rate; we compare it with that of 
Rosenblatt's estimator $\tilde f_n$ for which laws of the iterated logarithm 
were established by Hall (1981) in the case $d=1$ and by 
Arcones (1997) in the multivariate framework. \\

The remainder of the paper is organized as follows. In Section  \ref{section 2}, we state our main results: the bias and variance of $f_n$ are given in Subsection \ref{Section 2.0}, the 
pointwise estimation is considered in Subsection \ref{section 2.1}, the estimation by confidence intervals is developed in Subsection \ref{section 2.2}, and the strong convergence rate of $f_n$ is stated in Subsection \ref{Section 2.3}. 
Section \ref{section simul} is devoted to our simulation results, and 
Section \ref{section 5}  to the proof of our theoritical results. 

\section{Assumptions and main results} \label{section 2}

We consider stepsizes and bandwidths, which belong to the following 
class of regularly varying sequences.
\begin{defi}
Let $\gamma \in \mathbb{R} $ and $\left(v_n\right)_{n\geq 1}$ be a nonrandom positive sequence. We say that $\left(v_n\right) \in \mathcal{GS}\left(\gamma \right)$ if
\begin{eqnarray}\label{eq:5}
\lim_{n \to +\infty} n\left[1-\frac{v_{n-1}}{v_{n}}\right]=\gamma .
\end{eqnarray}
\end{defi}
Condition~\eqref{eq:5} was introduced by Galambos and Seneta (1973) to define regularly varying sequences (see also Bojanic and Seneta (1973)), and by Mokkadem and Pelletier (2007a) in the context of stochastic approximation algorithms. Typical sequences in $\mathcal{GS}\left(\gamma \right)$ are, for $b\in \mathbb{R}$, $n^{\gamma}\left(\log n\right)^{b}$, $n^{\gamma}\left(\log \log n\right)^{b}$, and so on. \\

The assumptions to which we shall refer are the following.
\begin{description}
\item(A1) $K:{\mathbb R^d}\rightarrow {\mathbb R}$ is a continuous, bounded 
function satisfying $\int_{\mathbb{R}^d}K\left( z\right) dz=1$, and, for 
all $j\in\{1,\ldots d\}$, $\int_{\mathbb{R}}z_jK\left( z\right) dz_j=0$ 
and $\int_{\mathbb{R}^d}z_j^2|K\left( z\right)| dz<\infty$. 
\item(A2) $i)$  $\left(\gamma_n\right) \in \mathcal{GS}\left(-\alpha \right)$ with 
$\alpha\in \left]{1}/{2},1\right]$.  \\
  $ii)$ $\left( h_{n}\right)\in \mathcal{GS} \left(-a\right)$ with 
$a \in \left]0,{\alpha}/{d}\right[$.\\ 
 $iii)$  $\lim_{n\to\infty}\left(n\gamma_n\right)\in]\min\{2a,(1-ad)/2\},\infty]$.
\item(A3) $f$ is bounded, twice differentiable, and, for all $i,j\in\{1,\ldots d\}$, 
$\partial^2f/\partial x_i\partial x_j$ is bounded.
\end{description}
Assumption $\left( A2\right)iii)$ on the limit of $\left(n\gamma_n\right)$ as $n$ goes to infinity is usual in the framework of stochastic approximation algorithms. It implies in particular that 
the limit of $\left([n\gamma_n]^{-1}\right)$ is finite.
Throughout this paper we will use the following notation:
\begin{eqnarray}
&&\xi= \lim_{n\to +\infty}\left(n\gamma_n\right)^{-1},\label{eq:6}\\
&&\mu_j^2=\int_{\mathbb{R}^d}z_j^2K\left(z\right)dz,\nonumber \\
&&f^{(2)}_{ij}(x) = \frac{\partial^2f}{\partial x_i\partial x_j}(x).\nonumber
\end{eqnarray}

\subsection{Bias and Variance} \label{Section 2.0}

Our first result is the following proposition, which gives the bias and the variance of $f_n$.

\begin{prop}[{Bias and Variance of \boldmath${f_n}$}]\label{Pr:1}
Let Assumptions $\left( A1\right)-\left( A3\right)$ hold, and assume that, 
for all $i,j\in\{1,\ldots d\}$, $f^{\left(2\right)}_{ij}$ is continuous at $x$.
\begin{enumerate}
\item If $a\leq\alpha/(d+4)$, then
\begin{eqnarray}\label{eq:7}
\mathbb{E}\left( f_{n}\left( x\right) \right) -f\left( x\right)=\frac{1}{2\left(1-2a\xi\right)}h_n^2
\sum_{j=1}^d\left(\mu_j^2f^{(2)}_{jj}(x)\right)
+o\left(h_n^2\right).
\end{eqnarray}
If $a>\alpha/(d+4)$, then
\begin{eqnarray}\label{eq:8}
\mathbb{E}\left( f_{n}\left( x\right) \right) -f\left( x\right)=o\left(\sqrt{\gamma_nh_n^{-d}}\right).
\end{eqnarray}
\item If $a\geq\alpha/(d+4)$, then
\begin{eqnarray}\label{eq:9}
Var\left( f_{n}\left( x\right) \right) =\frac{1}{2-\left(1-ad\right)\xi}\frac{\gamma_n}{h_n^d}f\left(x\right)\int_{\mathbb{R}^d} K^2\left(z\right)dz+o\left(\frac{\gamma_n}{h_n^d}\right).
\end{eqnarray}
If $a<\alpha/(d+4)$, then
\begin{eqnarray}\label{eq:10}
Var\left( f_{n}\left( x\right) \right) =o\left(h_n^{4}\right).
\end{eqnarray}
\item If $\lim_{n\to\infty}\left(n\gamma_n\right)>\max\{2a,(1-ad)/2\}$, 
then~\eqref{eq:7} and~\eqref{eq:9} hold simultaneously.
\end{enumerate}
\end{prop}


The bias and the variance of the estimator $f_n$ defined by the stochastic 
approximation algorithm (\ref{eq:1}) thus heavily depends on the choice of the stepsize 
$(\g_n)$. Let us recall that the bias and variance of Rosenblatt's estimator 
$\tilde f_{n}$ are given by:
\begin{eqnarray}
&&\mathbb{E}\left( \tilde f_{n}\left( x\right) \right) -f\left( x\right)  =  
\frac{1}{2}h_n^2
\sum_{j=1}^d\left(\mu_j^2f^{(2)}_{jj}(x)\right)
+o\left(h_n^2\right), \label{biais Ros}  \\
&&Var\left( \tilde f_{n}\left( x\right) \right)  =  
\frac{1}{nh_n^d}f\left(x\right)\int_{\mathbb{R}^d} K^2\left(z\right)dz
+o\left(\frac{1}{nh_n^d}\right). \label{variance Ros}
\end{eqnarray}
To illustrate the results given by Proposition \ref{Pr:1}, we now give some 
examples of possible choices of $(\g_n)$, and compare the bias and variance of 
$f_n$ with those of $\tilde f_{n}$.

\paragraph{Example 1: Choices of $(\gamma_n)$ minimizing the bias of $f_n$}

In view of~\eqref{eq:7}, the asymptotic bias of $f_{n}\left( x\right)$ is minimum when 
$\xi=0$, that is, when $\left(\gamma_n\right)$ is chosen such that 
$\lim_{n\to\infty}\left(n\gamma_n\right)=\infty$, 
and  we then have
\begin{eqnarray*}
\bE\left[f_n\left(x\right)\right]-f(x)=
\frac{1}{2}h_n^2
\sum_{j=1}^d\left(\mu_j^2f^{(2)}_{jj}(x)\right)
+o\left(h_n^2\right).
\end{eqnarray*}
In view of~(\ref{biais Ros}),  
the order of the bias of the recursive estimator $f_n$ is thus always greater 
or equal to that 
of Rosenblatt's estimator. Let us also mention that choosing the stepsize 
such that $\lim_{n\tinf}n\g_n=\infty$ (in which case the bias of $f_n$ is 
equivalent to that of Rosenblatt's estimator) is absolutely unadvised since 
we then have 
$$\lim_{n\to\infty}\frac{Var\left( \tilde f_{n}\left( x\right) \right)}
{Var\left(f_{n}\left( x\right) \right)}=0.$$

\paragraph{Example 2: Choices of $(\gamma_n)$ minimizing the variance of $f_n$}

As mentioned in the introduction, it is advised to minimize the variance of $f_n$ 
for interval estimation. 

\begin{coro} \label{corollaire 3}
Let the assumptions of Proposition \ref{Pr:1} hold with
$f(x)>0$. 
To minimize the asymptotic variance of $f_n$, $\alpha$ must be chosen equal to $1$,  
$\left(\gamma_n\right)$ must satisfy $\lim_{n\to \infty}n\gamma_n=1-ad$, 
and we then have
\begin{eqnarray*}
Var\left[f_n\left(x\right)\right]=\frac{1-ad}{nh_n^d}f\left(x\right)
\int_{\mathbb{R}^d}K^2\left(z\right)dz+o\left(\frac{1}{nh_n^d}\right).
\end{eqnarray*}
\end{coro}

It follows from Corollary \ref{corollaire 3} and (\ref{variance Ros}) that,
thanks to an adequate choice of $(\g_n)$,  
the variance of the recursive estimator $f_n$ can be smaller than that 
of Rosenblatt's estimator. To see better the comparison with Rosenblatt's estimator, 
let us set 
$(h_n)\in{\cal GS}(-1/[d+4])$ (which is the choice leading in particular to the 
minimum mean squared error of Rosenblatt's estimator). When $(\g_n)$ is chosen 
in ${\cal GS}(-1)$ and such that 
$\lim_{n\to \infty}n\gamma_n=1-d/[d+4]$, we have 
\begin{eqnarray}
\label{exemple2}
\lim_{n\tinf}\frac{\mathbb{E}\left( \tilde f_{n}\left( x\right) \right) -f\left( x\right)}
{\mathbb{E}\left(f_{n}\left( x\right) \right) -f\left( x\right)}
=\frac{1}{2}, & & 
\lim_{n\tinf}\frac{Var\left( \tilde f_{n}\left( x\right) \right)}
{Var\left(f_{n}\left( x\right) \right)}=\frac{d+4}{4}.
\end{eqnarray}
It is interesting to note that, whatever the dimension $d$ is, the bias of the 
recursive estimator $f_n$ is equivalent to twice that of Rosenblatt's 
estimator, whereas the 
ratio of the variances goes to infinity as the dimension $d$ increases. 

To conclude this example, let us mention that the most simple stepsize satisfying 
the conditions required in Corollary \ref{corollaire 3} is $(\g_n)=([1-ad]n^{-1})$. 

\paragraph{Example 3: The class of recursive estimators introduced by 
Hall and Patil (1994)}

The following lemma 
ensures that 
Proposition~\ref{Pr:1} gives the bias and variance of the recursive estimators 
defined in \eqref{peter} and introduced by Hall and Patil (1994) for 
a large choice of weights $(w_n)$.

\begin{lemma}\label{L:0}
Set $(w_n)\in{\cal GS}(w^*)$ and $(\g_n)=(w_n[\sum_{k=1}^nw_k]^{-1})$. If 
$w^*>-1$, then $(\g_n)\in{\cal GS}(-1)$ 
and $\lim_{n\tinf}n\g_n=1+w^*$.
\end{lemma}

Set $(h_n)\in{\cal GS}(-a)$; we give explicitly here the bias and variance of three 
particular recursive estimators.

\begin{itemize}
\item 
When $(w_n)=1$, $f_n$ is the estimator introduced by Wolverton and Wagner (1969); 
in view of Lemma \ref{L:0}, Proposition~\ref{Pr:1} applies with $\xi=1$, and we 
have
\begin{eqnarray*}
& &\mathbb{E}\left( f_{n}\left( x\right) \right) -f\left( x\right)=\frac{1}
{2\left(1-2a\right)}h_n^2
\sum_{j=1}^d\left(\mu_j^2f^{(2)}_{jj}(x)\right)
+o\left(h_n^2\right),  \\
& &Var\left( f_{n}\left( x\right) \right) =\frac{1}{1+ad}\frac{1}{nh_n^d}f\left(x\right)\int_{\mathbb{R}^d} K^2\left(z\right)dz+o\left(\frac{\gamma_n}{h_n^d}\right).
\end{eqnarray*}
\item When $(w_n)=(h_n^{d/2})$, $f_n$ 
is the estimator considered by Wegman and Davies (1979); 
in view of Lemma \ref{L:0}, Proposition~\ref{Pr:1} applies with $\xi=(1-ad/2)^{-1}$, 
and we have 
\begin{eqnarray*}
& &\mathbb{E}\left( f_{n}\left( x\right) \right) -f\left( x\right)=\frac{2-ad}
{2\left(2-[4+d]a\right)}h_n^2
\sum_{j=1}^d\left(\mu_j^2f^{(2)}_{jj}(x)\right)
+o\left(h_n^2\right),  \\
& &Var\left( f_{n}\left( x\right) \right) =\frac{(2-ad)^2}{4}\frac{1}{nh_n^d}f\left(x\right)\int_{\mathbb{R}^d} K^2\left(z\right)dz+o\left(\frac{\gamma_n}{h_n^d}\right).
\end{eqnarray*}
\item When $(w_n)=(h_n^{d})$, $f_n$ is the estimator introduced by Deheuvels (1973) 
and whose convergence rate was established by Duflo (1997); in view of 
Lemma \ref{L:0}, 
Proposition~\ref{Pr:1} applies with $\xi=(1-ad)^{-1}$, and we have 
\begin{eqnarray*}
& &\mathbb{E}\left( f_{n}\left( x\right) \right) -f\left( x\right)=\frac{1-ad}
{2\left(1-[2+d]a\right)}h_n^2
\sum_{j=1}^d\left(\mu_j^2f^{(2)}_{jj}(x)\right)
+o\left(h_n^2\right),  \\
& &Var\left( f_{n}\left( x\right) \right) =\frac{1-ad}{nh_n^d}f\left(x\right)\int_{\mathbb{R}^d} K^2\left(z\right)dz+o\left(\frac{\gamma_n}{h_n^d}\right).
\end{eqnarray*}
Let us underline that the bias and variance of this estimator are equivalent to 
those of the estimator defined with the stepsize $(\g_n=([1-ad]n^{-1})$ (this choice 
minimizing the variance of $f_n$, \emph{see} Corollary \ref{corollaire 3}), 
but its updating is less straightforward. 
\end{itemize}

\subsection{Choice of the optimal stepsize for point estimation} \label{section 2.1}

We first explicit the choices of $(\gamma_n)$ and $(h_n)$, which minimize the MSE and 
MISE of the recursive estimator defined by the stochastic approximation 
algorithm~\eqref{eq:1}, and then provide a comparison with Rosenblatt's estimator.

\subsubsection{Choices of $(\gamma_n)$ minimizing the MSE of $f_n$}

\begin{coro}\label{C:1} 
Let Assumptions $\left( A1\right)-\left( A3\right)$ hold, assume that $f(x)>0$, 
$\sum_{j=1}^d\left(\mu_j^2f^{(2)}_{jj}(x)\right)\neq 0$, 
and that, for all $i,j\in\{1,\ldots d\}$,  
$f^{\left(2\right)}_{ij}$ is continuous at $x$. To minimize the MSE of $f_n$ at the point $x$, the stepsize $\left(\gamma_n\right)$ must be chosen in $\mathcal{GS}\left(-1\right)$ and such that $\lim_{n\to \infty}n\gamma_n=1$, the bandwidth $\left(h_n\right)$ must equal
\begin{eqnarray*}
\left(\left[\frac{d(d+2)}{2(d+4)}\frac{f\left(x\right)\int_{\mathbb{R}^d} K^2\left(z\right)dz}
{\left(\sum_{j=1}^d\mu_j^2f^{(2)}_{jj}(x)\right)^2}\right]^{\frac{1}{d+4}}\gamma_n^{\frac{1}{d+4}}\right),
\end{eqnarray*}
and we then have
\begin{eqnarray*}
MSE=n^{-\frac{4}{d+4}}\frac{(d+4)^{\frac{3d+8}{d+4}}}
{d^{\frac{d}{d+4}}4^{\frac{d+6}{d+4}}(d+2)^{\frac{2d+4}{d+4}}}
\left[\sum_{j=1}^d\mu_j^2f^{(2)}_{jj}(x)\right]^{\frac{2d}{d+4}}
\left[f\left(x\right)\int_{\mathbb{R}^d} K^2\left(z\right)dz\right]^{\frac{4}{d+4}}
\left[1+o\left(1\right)\right].
\end{eqnarray*}
\end{coro}

The most simple example of stepsize belonging to $\mathcal{GS}\left(-1\right)$ and such 
that $\lim_{n\to \infty}n\gamma_n=1$ is $\left(\gamma_n\right)=\left(n^{-1}\right)$. 
For this choice of stepsize, the estimator $f_n$ defined by~\eqref{eq:1} equals 
the recursive kernel estimator introduced by Wolverton and Wagner (1969). This lattest 
estimator thus belongs to the subclass of recursive kernel estimators, which, 
thanks to an adequate choice of the bandwidth, have a minimum MSE. \\

\subsubsection{Choices of $(\gamma_n)$ minimizing the MISE of $f_n$}

The following proposition gives the MISE of the estimator $f_n$.
\begin{prop} \label{Pr:2} 
Let Assumptions $\left( A1\right)-\left( A3\right)$ hold, and assume that, 
for all $i,j\in\{1,\ldots d\}$, $f^{\left(2\right)}_{ij}$ 
is continuous and integrable.
\begin{enumerate}
\item If $a<\alpha/(d+4)$, then
\begin{eqnarray*}
MISE=\frac{1}{4\left(1-2a\xi\right)^2}h_n^4\int_{\mathbb{R}^d}
\left[\sum_{j=1}^d\mu_j^2f^{(2)}_{jj}(x)\right]^2dx +o\left(h_n^4\right).
\end{eqnarray*}
\item If $a=\alpha/(d+4)$, then
\begin{eqnarray*}
MISE & = & \frac{1}{4\left(1-2a\xi\right)^2}h_n^4\int_{\mathbb{R}^d}
\left[\sum_{j=1}^d\mu_j^2f^{(2)}_{jj}(x)\right]^2dx 
+\frac{1}{2-\left(1-ad\right)\xi}\frac{\gamma_n}{h_n^d}\int_{\mathbb{R}^d} 
K^2\left(z\right)dz\\
& & \mbox{ }+o\left(h_n^4+\frac{\gamma_n}{h_n^d}\right).
\end{eqnarray*}
\item If $a>\alpha/(d+4)$, then
\begin{eqnarray*}
MISE=\frac{1}{2-\left(1-ad\right)\xi}\frac{\gamma_n}{h_n^d}\int_{\mathbb{R}^d} 
K^2\left(z\right)dz+o\left(\frac{\gamma_n}{h_n^d}\right).
\end{eqnarray*}
\end{enumerate}
\end{prop}

The following corollary ensures that Wolwerton and Wagner's estimator also belongs to the subclass 
of kernel estimators defined by the stochastic approximation algorithm~\eqref{eq:1}, which, thanks to 
an adequate choice of the bandwidth, have a minimum MISE.  

\begin{coro} \label{C:2} Let Assumptions $\left( A1\right)-\left( A3\right)$ hold, and assume 
that, for all $i,j\in\{1,\ldots d\}$, $f^{\left(2\right)}_{ij}$ is continuous 
and integrable. To minimize the MISE of $f_n$, the stepsize $\left(\gamma_n\right)$ must be chosen in $\mathcal{GS}\left(-1\right)$ and such that $\lim_{n\to \infty}n\gamma_n=1$, the bandwidth $\left(h_n\right)$ must equal
\begin{eqnarray*}
\left(\left[\frac{d(d+2)}{2(d+4)}\frac{\int_{\mathbb{R}^d} K^2\left(z\right)dz}
{\int_{\mathbb{R}^d}\left(\sum_{j=1}^d\mu_j^2f^{(2)}_{jj}(x)\right)^2dx}
\right]^{\frac{1}{d+4}}\gamma_n^{\frac{1}{d+4}}\right),
\end{eqnarray*}
and we then have 
\begin{eqnarray*}
MISE=n^{-\frac{4}{d+4}}\frac{(d+4)^{\frac{3d+8}{d+4}}}
{d^{\frac{d}{d+4}}4^{\frac{d+6}{d+4}}(d+2)^{\frac{2d+4}{d+4}}}
\left[\int_{\mathbb{R}^d}\left(\sum_{j=1}^d\mu_j^2f^{(2)}_{jj}(x)\right)^2dx
\right]^{\frac{d}{d+4}}\left[\int_{\mathbb{R}^d} K^2\left(z\right)dz\right]^{\frac{4}{d+4}}\left[1+o\left(1\right)\right].
\end{eqnarray*}
\end{coro}

\subsubsection{Comparison with Rosenblatt's estimator}

The ratio of the optimal MSE (or MISE) of Rosenblatt's estimator to that of 
Wolwerton and Wagner's estimator equals 
$$\rho(d)={\left[\frac{2^4(d+2)^{2d+4}}{(d+4)^{2d+4}}\right]}^{\frac{1}{d+4}}.$$
This ratio is always less than one, it at first decreases, and then increases to 
one as the dimension $d$ increases. This phenomenon is similar to that  
observed by Hall and Patil (1994). The former authors consider the univariate 
framework, but look at the efficiency of Wolwerton and Wagner's estimator 
of the $s$th-order derivative of $f$ relative to Rosenblatt's one; the ratio 
$\rho(s)$ varies in $s$ in the same way as $\rho(d)$ does in $d$. According to 
pointwise estimation point of view, and when rapid updating is not too important, 
it is thus preferable to use  Rosenblatt's nonrecursive estimator rather than 
any recursive estimator defined by the stochastic approximation 
algorithm~\eqref{eq:1}.  Let us mention that Hall and Patil (1994) introduce 
a class of on-line estimators, constructed from the class of 
the recursive estimators defined in \eqref{peter}; their on-line estimators are not 
recursive any more, but updating them requires much less operations than updating 
Rosenblatt's estimator, and their MSE and MISE are smaller than those of the 
recursive estimators \eqref{peter}.

\subsection{Choice of the optimal stepsize for interval estimation} \label{section 2.2}

Let us first state the following theorem, which gives the weak convergence rate of the estimator $f_n$ defined in~\eqref{eq:1}.

\begin{theor}[Weak pointwise convergence rate]\label{T:1} 
Let Assumptions $\left( A1\right)-\left( A3\right) $ hold, assume that $f(x)>0$ 
and that, for all $i,j\in\{1,\ldots d\}$, $f^{\left(2\right)}_{ij}$ is 
continuous at $x$.
\begin{enumerate}
\item If there exists $c\geq 0 $ such that $\gamma_n^{-1}h_{n}^{d+4}\rightarrow c$, then
\begin{eqnarray*}
\lefteqn{\sqrt{\gamma_n^{-1}h_{n}^d}\left( f_{n}\left( x\right)-f\left( x\right) \right)}
\\ & \stackrel{\mathcal{D}}{\rightarrow} &
\mathcal{N}\left( \frac{c^{\frac{1}{2}}}{2\left(1-2a\xi\right)}
\sum_{j=1}^d\left(\mu_j^2f^{(2)}_{jj}(x)\right)
,\frac{1}{2-\left(1-ad\right)\xi}f\left(x\right) \int_{\mathbb{R}^d}K^2\left(z\right)dz\right).
\end{eqnarray*}
\item If $\gamma_n^{-1}h_{n}^{d+4} \rightarrow \infty $, then  
\begin{eqnarray*}
\frac{1}{h_{n}^{2}}\left( f_{n}\left( x\right)-f\left( x\right) \right) \stackrel{\mathbb{P}}{\rightarrow } \frac{1}{2\left(1-2a\xi\right)}\sum_{j=1}^d\left(\mu_j^2f^{(2)}_{jj}(x)\right),
\end{eqnarray*}
\end{enumerate}
where $\stackrel{\mathcal{D}}{\rightarrow}$ denotes the convergence in distribution, $\mathcal{N}$ the Gaussian-distribution and $\stackrel{\mathbb{P}}{\rightarrow}$ the convergence in probability.
\end{theor}

As mentioned in the introduction, Hall (1992) shows that, to minimize the coverage error of probability density confidence intervals, avoiding bias estimation by a slight undersmoothing is more efficient than bias correction. 
Let us recall that, when the bandwidth $(h_n)$ is chosen such that 
$\lim_{n\to\infty}nh_n^{d+4}=0$ (which corresponds to undersmoothing), 
Rosenblatt's estimator fulfills the central limit theorem 
\begin{eqnarray}
\label{tlc}
\sqrt{nh_{n}^d}\left( \tilde f_{n}\left( x\right)-f\left( x\right) \right) 
& \stackrel{\mathcal{D}}{\rightarrow} & \mathcal{N}\left(0,f\left(x\right) 
\int_{\mathbb{R}^d}K^2\left(z\right)dz\right).
\end{eqnarray}
Now, let $\Phi$ denote the distribution function of the $\mathcal{N}(0,1)$, 
let $t_{\alpha/2}$ be such that $\Phi(t_{\alpha/2})=1-\alpha/2$ (where 
$\alpha\in]0,1[$), and set 
\begin{eqnarray*}
I_{g_n}(x)&=&\left[g_n\left(x\right)-t_{\alpha/2}C\left(g_n\right)
\sqrt{\frac{g_n\left(x\right)\int_{\mathbb{R}^d}K^2\left(z\right)dz}{nh_n^d}},\,g_n\left(x\right)
+t_{\alpha/2}C\left(g_n\right)
\sqrt{\frac{g_n\left(x\right)\int_{\mathbb{R}^d}K^2\left(z\right)dz}{nh_n^d}}\right].
\end{eqnarray*}
In view of (\ref{tlc}), the asymptotic level of $I_{\tilde f_n}(x)$ equals $1-\alpha$ 
for $C(\tilde f_n)=1$. The following corollary gives the values of $C(f_n)$ for which 
the asymptotic level of $I_{f_n}(x)$ equals $1-\alpha$ too.

\begin{coro} \label{corollaire niveau}
Let the assumptions of Theorem \ref{T:1} hold with 
$\lim_{n\to \infty}n\gamma_n=\g_0\in]0,\infty[$ and $\lim_{n\to \infty}nh_n^{d+4}=0$. 
The asymptotic level of $I_{f_n}(x)$ equals $1-\alpha$ for $$C(f_n)=\sqrt{\gamma_0
\left[2-\left(1-ad\right)\gamma_0^{-1}\right]^{-1}}.$$
Moreover, the minimum of $C(f_n)$ is reached at $\g_0=1-ad$ and equals $\sqrt{1-ad}$. 
\end{coro}

The optimal stepsizes for interval estimation are thus the sequences $(\g_n)\in{\cal GS}(-1)$ 
such that  $\lim_{n\to \infty}n\gamma_n=1-ad$, the most simple one being $(\g_n)=([1-ad]n^{-1})$. 
Of course, these stepsizes are those which minimize the variance of $f_n$ (see Corollary 
\ref{corollaire 3}). 


\subsection{Strong pointwise convergence rate} \label{Section 2.3}

The following theorem gives the strong pointwise convergence rate of $f_n$.

\begin{theor}[Strong pointwise convergence rate]\label{T:2} 
Let Assumptions $\left( A1\right)-\left( A3\right) $ hold, and assume that, 
for all $i,j\in\{1,\ldots d\}$, $f^{\left(2\right)}_{ij}$ is continuous at $x$.
\begin{enumerate}
\item If there exists $c_{1}\geq 0 $ such that $\gamma_n^{-1}h_{n}^{d+4}/\left(\ln[\sum_{k=1}^n\g_k]\right)\rightarrow c_{1}$, then, with probability one, the sequence
\[
\left(\sqrt{\frac{\gamma_n^{-1}h_{n}^d}{2\ln [\sum_{k=1}^n\g_k]}}
\left( f_{n}\left( x\right)-f\left( x\right) \right) \right) 
\]
is relatively compact and its limit set is the interval
\begin{eqnarray*}
\lefteqn{\left[\frac{1}{2\left(1-2a\xi\right)}\sqrt{\frac{c_1}{2}}
\sum_{j=1}^d\left(\mu_j^2f^{(2)}_{jj}(x)\right)
-\sqrt{\frac{f\left(x\right)}{2-\left(1-ad\right)\xi}
\int_{\mathbb{R}^d} K^2\left(z\right)dz},\right.}\\
&&\left.\frac{1}{2\left(1-2a\xi\right)}\sqrt{\frac{c_1}{2}}
\sum_{j=1}^d\left(\mu_j^2f^{(2)}_{jj}(x)\right)
+\sqrt{\frac{f\left(x\right)}{2-\left(1-ad\right)\xi}
\int_{\mathbb{R}^d} K^2\left(z\right)dz}\right].
\end{eqnarray*}

\item If $\gamma_n^{-1}h_{n}^{d+4}/\left( \ln [\sum_{k=1}^n\g_k]\right) 
\rightarrow \infty $, then, with probability one,
\begin{eqnarray*}
\lim_{n\rightarrow \infty}\frac{1}{h_{n}^{2}}\left( f_{n}\left( x\right)-f\left( x\right) \right) = \frac{1}{2\left(1-2a\xi\right)}
\sum_{j=1}^d\left(\mu_j^2f^{(2)}_{jj}(x)\right).
\end{eqnarray*}
\end{enumerate}
\end{theor}

Set $(h_n)$ such that $\lim_{n\tinf}nh_{n}^{d+4}/\ln\ln n=0$. Arcones (1997) proves the 
following compact law of the iterated logarithm for Rosenblatt's estimator:
with probability one, the sequence 
$(\sqrt{nh_{n}^d}( \tilde f_{n}( x)-f( x))/\sqrt{2\ln\ln n})$ is relatively compact 
and its limit set is the interval 
$$J=\left[
-\sqrt{f\left(x\right)\int_{\mathbb{R}^d} K^2\left(z\right)dz},
\sqrt{f\left(x\right)\int_{\mathbb{R}^d} K^2\left(z\right)dz}\right].$$ 
Now, set $(\g_n)$ such that $\lim_{n\tinf}n\g_{n}=\g_0\in]0,\infty[$. 
The first part of Theorem \ref{T:2} ensures that, with probability one, 
the limit set of the sequence 
$(\sqrt{nh_{n}^d}( f_{n}( x)-f( x))/\sqrt{2\ln\ln n}) $ 
is the interval 
$$J(\g_0)=\left[
-A(\g_0)\sqrt{f\left(x\right)\int_{\mathbb{R}^d} K^2\left(z\right)dz},
A(\g_0)\sqrt{f\left(x\right)\int_{\mathbb{R}^d} K^2\left(z\right)dz}\right]
\mbox{ with }
A(\g_0)=\sqrt{\frac{\g_0}{[2-(1-ad)\g_0^{-1}]}}.$$ 
In particular, for Wolwerton and Wagner's estimator, $A(\g_0)=1/\sqrt{1+ad}$; 
for the estimator considered by Wegman and Davies (1979), 
or when $(\g_n)=([1-ad/2]n^{-1})$,  $A(\g_0)=1-ad/2$; 
for the estimator considered by Deheuvels (1973) and Duflo (1997), 
or when $(\g_n)=([1-ad]n^{-1})$,  $A(\g_0)=\sqrt{1-ad}$. 
For all these recusive estimators, the length of the limit interval $J(\g_0)$ is smaller 
than that of $J$, which shows that they are more concentrated around $f$ than 
Rosenblatt's estimator is.

\section{Simulations} \label{section simul}

The aim of our simulation studies is to compare the performance of Rosenblatt's 
estimator defined in~\eqref{eq:4} with that of the recursive estimators, from 
confidence interval  point of view. Of course, the recursive estimator we consider 
here is the optimal one according to this criteria (see Corollary 
\ref{corollaire niveau}). We set: 
\begin{eqnarray*}
I_{i,n}&=&\left[g_n\left(x\right)-1.96\ C(g_n)
\sqrt{\frac{g_n\left(x\right)\int_{\mathbb{R}^d}K^2\left(z\right)dz}{nh_n^d}},\,g_n\left(x\right)+1.96\ C(g_n)\sqrt{\frac{g_n\left(x\right)\int_{\mathbb{R}^d}K^2\left(z\right)dz}{nh_n^d}}\right],
\end{eqnarray*}
where:
\begin{itemize}
\item if $i=1$, then $g_n=\ti f_n$ is Rosenblatt's estimator, and $C(g_n)=1$;
\item if $i=2$, then $g_n=f_n$ is the optimal recursive estimator defined by 
the algorithm~\eqref{eq:1} with the stepsize $\left(\gamma_n\right)
=\left([1-ad]n^{-1}\right)$, and $C(g_n)=\sqrt{1-ad}$.
\end{itemize}
According to the theoritical results given in Section \ref{section 2.2}, both 
confidence intervals $I_{1,n}$ and $I_{2,n}$ have the same asymptotic level 
(equal to 95\%), whereas $I_{2,n}$ has a smaller length than $I_{1,n}$. 
In order to investigate their finite sample behaviours, we consider three 
sample sizes: $n=50$, $n=100$, and $n=200$.
In each case, the number of simulations is $N=5000$.  
Tables 1-4 give (for different values of $d$, $f$, $x$, and $(h_n)$): 
\begin{itemize}
\item 
the empirical levels $\#\left\{f\left(x\right)\in I_{i,n}\right\}/N$ at each 
first line concerning $I_{i,n}$. 
\item
the averaged lengths of the intervals $I_{i,n}$ at each 
second line concerning $I_{i,n}$. 
\end{itemize}

\paragraph{The case \boldmath $d=1$.}
In the univariate framework, we consider two densities $f$: the standard normal 
$\mathcal{N}(0,1)$ distribution (see Table 1), and the normal mixture 
$\frac{1}{2}\mathcal{N}(-\frac{1}{2},1)+\frac{1}{2}\mathcal{N}(\frac{1}{2},1)$ 
distribution (see Table 2). 
The points at which $f$ is estimated are: $x=0$, $x=0.5$, and $x=1$. 
The bandwidth $(h_n)$ is set equal to $(n^{-a})$ with $a=0.21$ and $a=0.23$ (the 
parameter $a$ being chosen slightly larger than $1/5$ to slightly undersmooth).
Both tables show that the recursive estimator performs better than 
Rosenblatt's one: the empirical levels of the intervals $I_{2,n}$ are greater 
than those of $I_{1,n}$, whereas their averaged lengths are smaller. 

\paragraph{The case \boldmath $d=2$.}
In the case when $d=2$, we estimate the density $f$ of the random vector $X$ defined 
as $X=AY$ with 
$A=\begin{pmatrix}
1 &   0\\
0.5  & 1
\end{pmatrix}$, and where the distribution of the random vector $Y$ is:
\begin{itemize}
\item the normal standard distribution $\mathcal{N}\left(0,{I}_2\right)$ 
(see Table 3); 
\item the normal mixture $\frac{1}{2}\mathcal{N}\left(-B,{I}_{2}\right)
+\frac{1}{2}\mathcal{N}\left(B,{I}_{2}\right)$
with $B=\begin{pmatrix}
-0.5 \\
-0.5 
\end{pmatrix}$ (see Table 4).
\end{itemize}
The points at which $f$ is estimated are: $x=\left(0,0\right)$, 
$x=\left(0.5,0.5\right)$, and $x=\left(1,1\right)$. 
The bandwidth $(h_n)$ is set equal to $(n^{-a})$. To slightly undersmooth, the 
parameter $a$ must be chosen slightly larger than $1/6$; we first chose 
$a=0.17$ and $a=0.19$. Tables 3 and 4 show that, for these given values of the parameter $a$, 
the recursive estimator performs better for the sample sizes 
$n=50$ and $n=100$, whereas, at first glance, Rosenblatt's estimator performs 
better in the case 
when $n=200$. This is explained by the fact that, for this lattest sample size, 
the length of $I_{2,n}$ becomes too small. We have thus added other choices of 
the parameter $a$ ($a=0.21$ in Table 3; $a=0.21$ and $a=0.24$ in Table 4). 
The larger $a$ is, the larger the length of the intervals $I_{i,n}$ are, and 
the larger the empirical levels are. Now, Tables 3 and 4 also show that, 
for the sample size $n=200$, 
the intervals $I_{2,n}$ computed with $a=0.21$ or $a=0.24$ have a 
smaller length and a higher level than the intervals $I_{1,n}$ computed 
with $a=0.17$ or $a=0.19$, so that we can say again that the recursive 
estimator performs better than Rosenblatt's one. \\

This simulation study shows the good performance of the recursive estimator defined by 
the algorithm~\eqref{eq:1} with the stepsize $\left(\gamma_n\right)
=\left([1-ad]n^{-1}\right)$ for interval estimation. The main question which 
remains open is how to choose the bandwidth $(h_n)$ in $\mathcal{GS}(-a)$, and, 
in particular, how to determine the parameter $a$. This problem is not particular to 
the framework of recursive estimation; in the case when Rosenblatt's estimator is 
used, Hall (1992) enlightens that criteria to determine the ``good undersmoothing'' 
are not easy to determine empirically.

\begin{table}
\caption{$X\leadsto\mathcal{N}(0,1)$}
\begin{eqnarray*}
\begin{tabular}{lccc|ccc|ccc}
\label{gtable1}
& \multicolumn{3}{c|}{$x=0$}&\multicolumn{3}{c|}{$x=0.5$}&\multicolumn{3}{c}{$x=1$}\\
& $n=50$  & $n=100$ & $n=200$ &$n=50$  & $n=100$ & $n=200$ &  $n=50$  & $n=100$ & 
$n=200$  \\ \hline
&&  &&&\boldmath{$a=0.21$}  &&&& \\
$I_{1,n}$ &  96.74\% &  96.08\%  & 95.74\% & 97.1\% &  96.74\% & 96.96\% &97.72\% & 97.44\% &  97.7\%\\
& 0.2681& 0.2061& 0.158&0.2538& 0.1948&0.1493& 0.2168& 0.165&0.126\\
$I_{2,n}$&  99.36\% &  98\% & 96.18\% &99.76\% &   98.96\% &   98.36\%&98.86\% & 98.76\%& 98.78\%\\
& 0.2436 & 0.184&0.140&0.2332& 0.1755&0.1331&0.2068&0.1529&0.1146\\
\hline
&&  &&& \boldmath{$a=0.23$}  &&&& \\
$I_{1,n}$ &96.58\% &    96.46\% &   96.78\% &96.78\% & 97.06\%  &97.04\%&97.32\%  &97.58\%&96.96\% \\
&0.2796& 0.2167&  0.1674&0.2653& 0.205&0.1579&0.225&0.1731&0.1328\\
$I_{2,n}$ &99.46\% &  98.58\% & 97.58\% &99.6\% &   99.26\% & 98.72\% &98.68\% &98.32\%& 97.96\%\\
&0.2517& 0.1915& 0.1467&0.2415&0.1828&0.1393&0.2134& 0.159&0.1197\\
\end{tabular}
\end{eqnarray*}
\end{table}

\begin{table}
\caption{$X\leadsto\frac{1}{2}\mathcal{N}(-\frac{1}{2},1)+\frac{1}{2}\mathcal{N}
(\frac{1}{2},1)$}
\begin{eqnarray*}
\begin{tabular}{lccc|ccc|ccc}
 &\multicolumn{3}{c|}{$x=0$}&\multicolumn{3}{c|}{$x=0.5$}&\multicolumn{3}{c}{$x=1$}\\
 &$n=50$  & $n=100$ & $n=200$ &$n=50$  & $n=100$ & $n=200$ &  $n=50$  & $n=100$ & $n=200$  \\ \hline
&&  &&&\boldmath{$a=0.21$}  &&&& \\
$I_{1,n}$&   96.86\% &  96.96\%  &96.86\% & 96.96\%&96.68\% & 96.8\% &97.12\% & 97.04\% &  96.94\%\\
& 0.2541&0.1949&0.1493&0.2436&0.1866&0.1427&0.2142& 0.1642&0.1251\\
$I_{2,n}$ & 99.76\% &  99.04\% & 98.2\% &99.62\% &   99.28\% &   98.72\%&99.14\% & 98.94\%& 98.4\%\\
& 0.2334 &0.1755&0.1331&0.2257&0.1692&0.1278&0.2045&0.1518&0.1136\\
\hline
&&  &&&\boldmath{$a=0.23$}  &&&& \\ 
$I_{1,n}$&96.92\% & 97.04\% & 96.84\% &96.56\% & 96.66\%  &97.14\%&97.02\%  &97.12\%&96.76\% \\
&0.2654&0.2049&0.1579&0.254&0.196&0.151&0.2233&0.1717&0.1321\\
$I_{2,n}$&99.9\% &  99.18\% & 98.76\% &99.74\% &   99.3\% & 98.92\% &98.78\% &98.76\%& 98.2\%\\
&0.2416&0.1826&0.1393&0.2334&0.176&0.1338&0.2116&0.1575&0.1187\\
\end{tabular}
\end{eqnarray*}
\end{table}

\begin{table}
\caption{$X=AY$ with $Y\leadsto\mathcal{N}\left(0,{I}_2\right)$}
\begin{eqnarray*}
\begin{tabular}{lccc|ccc|ccc}
& \multicolumn{3}{c|}{$x=\left(0,0\right)$}&\multicolumn{3}{c|}{$x=\left(0.5,0.5\right)$}&\multicolumn{3}{c}{$x=\left(1,1\right)$}\\
& $n=50$  & $n=100$ & $n=200$ &$n=50$  & $n=100$ & $n=200$ &  $n=50$  & $n=100$ & $n=200$  \\\hline
&&  &&&\boldmath{$a=0.17$}  &&&& \\
$I_{1,n}$&    93.82\% &  94.98\%  & \colorbox{light}{96.9\%} & 91.06\% & 92.82\% &
\colorbox{light}{94.0\%} &89.48\% & 86.88\% & 85.82\%\\
& 0.1159&0.0934&\colorbox{light}{0.0757}&0.1059&0.0854&
\colorbox{light}{0.0686}&0.0811&0.0645&0.0515\\
$I_{2,n}$&   97.54\% &  95.12\% & {94.34\%} &96.74\% &   94.62\% &   
{92.86\%}&97.2\% & 94.32\%& 91.16\%\\
&  0.0979&0.0765&{0.061}&0.091&0.0707&
{0.0558}&0.0736&0.0557&0.0432\\
\hline
&&  &&&\boldmath{$a=0.19$}  &&&& \\
$I_{1,n}$&   95.64\% &  97.08\%  & \colorbox{dark}{97.28\%} & 93.46\% & 94.84\% &95.82\% &91.58\% & 91.06\% &89.04\%\\
& 0.1271&0.1042&\colorbox{dark}{0.0851}&0.1158&0.0946&0.077&0.0883&0.0713&0.0574\\
$I_{2,n}$&   97.5\% &  97.26\% & {96.64\%} &97.22\% &96.5\% & 95.42\%&96.74\% & 95.66\%& 92.24\%\\
& 0.1045&0.0829&{0.0666}&0.0969&0.0763&0.0609&0.0783&0.0599&0.0469\\
\hline
&&  &&&\boldmath{$a=0.21$}  &&&& \\
$I_{1,n}$& 96.68\% & 97.62\% & 98.24\% &95.16\% & 96.48\%  &97.16\%&92.76\%  &91.2\%&91.04\% \\
&0.1392&0.1157&0.0957&0.1267&0.105&0.0863&0.0962&0.0783&0.0641\\
$I_{2,n}$& 97.16\% &  97.48\% & \colorbox{midmid}{97.56\%} &96.96\% &   96.84\% & 
\colorbox{light}{96.7\%} &96.72\% &96.58\%& 94.2\%\\
&0.1111&0.0893&\colorbox{midmid}{0.0726}&0.1031&0.0822&
\colorbox{light}{0.0662}&0.0832&0.0642&0.0509\\
\end{tabular}
\end{eqnarray*}
\end{table}

\begin{table}
\caption{$X=AY$ with $Y\leadsto\frac{1}{2}\mathcal{N}\left(-B,{I}_{2}\right)
+\frac{1}{2}\mathcal{N}\left(B,{I}_{2}\right)$}
\begin{eqnarray*}
\begin{tabular}{lccc|ccc|ccc}
& \multicolumn{3}{c|}{$x=\left(0,0\right)$}&\multicolumn{3}{c|}{$x=\left(0.5,0.5\right)$}&\multicolumn{3}{c}{$x=\left(1,1\right)$}\\
& $n=50$  & $n=100$ & $n=200$ &$n=50$  & $n=100$ & $n=200$ &  $n=50$  & $n=100$ & $n=200$  \\ \hline
&&  &&&\boldmath{$a=0.17$}  &&&& \\
$I_{1,n}$&     91.84\% &  91.28\%  & \colorbox{light}{92.4\%} & 90.06\% & 89.42\% &
\colorbox{light}{87.86\%} &83.24\% & 80.46\% & 78.88\%\\
& 0.105&0.0847&\colorbox{light}{0.068}&0.0976&0.0785&
\colorbox{light}{0.063}&0.0787&0.0631&0.050\\
$I_{2,n}$&    96.8\% &  93.76\% & {91.34\%} &95.9\% &   92.32\% &  
{86.96\%}&95.52\% & 87.6\%& 82.12\%\\
&  0.0903&0.0702&{0.0553}&0.0851&0.0657&
{0.0516}&0.0716&0.0544&0.0419\\
\hline
&&  &&&\boldmath{$a=0.19$}  &&&& \\
$I_{1,n}$&    93.54\% &  93.94\%  & \colorbox{dark}{95.44\%} & 90.72\% & 91.38\% &
\colorbox{dark}{92.12\%} &85.46\% & 84.24\% &82.24\%\\
& 0.1151&0.094&\colorbox{dark}{0.0764}&0.1158&0.1069&
\colorbox{dark}{0.0706}&0.0857&0.0692&0.0457\\
$I_{2,n}$&    97.42\% &  95.92\% & {94.38\%} &97.22\% &97.06\% & 
{91.74\%}&96.18\% & 91.26\%& 86.88\%\\
& 0.0964&0.0757&{0.0604}&0.0969&0.0908&
{0.0562}&0.0762&0.0582&0.0469\\
\hline
&&  &&&\boldmath{$a=0.21$}  &&&& \\
$I_{1,n}$&  94.82\% & 96.12\% & 97.44\% &93.14\% & 93.46\%  &94.16\%&88.72\%  &86.24\%&83.54\% \\
&0.1259&0.1037&0.0858&0.1163&0.0962&0.0793&0.0935&0.0764&0.0624\\
$I_{2,n}$&  97.1\% &  97.48\% & \colorbox{light}{96.96\%} &96.82\% &   96.04\% & 
\colorbox{light}{93.96\%} &96.76\% &93.52\%& 88.24\%\\
&0.1025&0.0813&\colorbox{light}{0.0659}&0.0963&0.0762&
\colorbox{light}{0.0613}&0.0811&0.0627&0.0495\\
\hline
&&  &&&\boldmath{$a=0.24$}  &&&& \\
$I_{1,n}$&  96.26\% & 97.48\% & 98.38\% &94.36\% & 96.16\%  &96.7\%&91.04\%  &91.08\%&89.42\% \\
&0.1435&0.1208&0.1017&0.1325&0.1117&0.0937&0.1058&0.0885&0.0736\\
$I_{2,n}$&   96.18\% &  97.54\% & \colorbox{dark}{98.04\%} &96.68\% & 97.38\% & 
\colorbox{dark}{96.6\%} &96.98\% &95.96\%& 91.3\%\\
&0.1117&0.0903&\colorbox{dark}{0.0743}&0.1049&0.0845&
\colorbox{dark}{0.0691}&0.0883&0.0695&0.0558\\
\end{tabular}
\end{eqnarray*}
\end{table}

\newpage

\section{Proofs} \label{section 5}
Throughout this section we use the following notation:
\begin{eqnarray}
&& \Pi_n=\prod_{j=1}^{n}\left(1-\gamma_j\right),\nonumber\\
&& s_n=\sum_{k=1}^n\gamma_k, \nonumber\\
&&Z_n\left(x\right)=\frac{1}{h_n^d}K\left(\frac{x-X_n}{h_n}\right).\label{eq:13}
\end{eqnarray}
Let us first state the following technical lemma. 

\begin{lemma}\label{L:1} 
Let $\left(v_n\right)\in \mathcal{GS}\left(v^*\right)$, 
$\left(\gamma_n\right)\in \mathcal{GS}\left(-\alpha\right)$, and $m>0$ 
such that $m-v^*\xi>0$ where $\xi$ is defined in~\eqref{eq:6}. We have 
\begin{eqnarray}\label{eq:14}
\lim_{n \to +\infty}v_n\Pi_n^{m}\sum_{k=1}^n\Pi_k^{-m}\frac{\gamma_k}{v_k}
=\frac{1}{m-v^*\xi}. 
\end{eqnarray}
Moreover, for all positive sequence $\left(\alpha_n\right)$ such that $\lim_{n \to +\infty}\alpha_n=0$, and all $\delta \in \mathbb{R}$,
\begin{eqnarray}\label{eq:15}
\lim_{n \to +\infty}v_n\Pi_n^{m}\left[\sum_{k=1}^n \Pi_k^{-m} \frac{\gamma_k}{v_k}\alpha_k+\delta\right]=0.
\end{eqnarray}
\end{lemma}

Lemma \ref{L:1} is widely applied throughout the proofs. Let us underline that it is 
its application, which requires Assumption $(A2)iii)$ on the limit of $(n\g_n)$ as $n$ goes to 
infinity. Let us mention that, in particular, to prove (\ref{eq:9}), Lemma \ref{L:1} 
is applied with $m=2$ and $(v_n)=(\g_n^{-1}h_n^d)$ (and thus $v^*=\alpha -ad$); 
the stepsize $(\g_n)$ must thus fulfill the condition 
$\lim_{n\to\infty}\left(n\gamma_n\right)>\left(\alpha -ad\right)/2$. Now, since 
$\lim_{n\to\infty}\left(n\gamma_n\right)<\infty$ only if $\alpha =1$, 
the condition $\lim_{n\to\infty}\left(n\gamma_n\right)\in]\min\{2a,(1-ad)/2\},\infty]$ in 
$(A2)iii)$ is equivalent to the condition 
$\lim_{n\to\infty}\left(n\gamma_n\right)\in]\min\{2a,(\alpha-ad)/2\},\infty]$, which appears 
throughout our proofs.
Similarly, since $\xi\neq 0$ only if $\alpha =1$, 
the limit $[2-\left(\alpha -ad\right)\xi]^{-1}$ given by the application of 
Lemma \ref{L:1} for such $m$ and $(v_n)$ equals the factor 
$[2-\left(1 -ad\right)\xi]^{-1}$ that stands in the statement of our main results.\\

Our proofs are now organized as follows. 
Lemmas \ref{L:0} and \ref{L:1} are proved in Section \ref{Section 5.0}, 
Propositions \ref{Pr:1} and \ref{Pr:2} in Sections \ref{Section 5.2} 
and \ref{Section 5.3} respectively, Theorems \ref{T:1} and \ref{T:2}
in Sections \ref{Section 5.4} and \ref{Section 5.5} respectively, 
and Corollaries \ref{corollaire 3}-\ref{corollaire niveau} in 
Section \ref{Section 5.2bis}.


\subsection{Proof of Lemmas~\ref{L:0} and \ref{L:1}} \label{Section 5.0}

We first prove Lemma~\ref{L:0}. 
Since $\left(w_n\right)\in\mathcal{GS}\left(w^*\right)$ with $w^*>-1$, we have
\begin{eqnarray}\label{eq:12}
\lim_{n\to\infty}\frac{nw_n}{\sum_{k=1}^nw_k}=1+w^*,
\end{eqnarray}
which guarantees that $\lim_{n\to\infty}n\gamma_n=1+w^*$. 
Moreover, applying~\eqref{eq:12}, we note that
$$
\frac{\sum_{k=1}^{n-1}w_k}{\sum_{k=1}^nw_k}=
1-\frac{w_n}{\sum_{k=1}^nw_k}
 =  1-\frac{1+w^*}{n}+o\left(\frac{1}{n}\right),
$$
so that
\begin{eqnarray*}
\lim_{n\to \infty}n\left[1-\frac{\sum_{k=1}^{n-1}w_k}{\sum_{k=1}^nw_k}\right]=1+w^*.
\end{eqnarray*}
It follows that $\left(\sum_{k=1}^nw_k\right)\in\mathcal{GS}\left(1+w^*\right)$, 
and thus that $\left(\gamma_n\right)\in \mathcal{GS}\left(-1\right)$, which 
concludes the proof of Lemma~\ref{L:0}.\\

To prove Lemma~\ref{L:1}, we first establish~\eqref{eq:15}. Set
\begin{eqnarray*}
Q_n=v_n\Pi_n^m\left[\sum_{k=1}^n\Pi_k^{-m}\gamma_kv_k^{-1}\alpha_k+\delta\right].
\end{eqnarray*}
We have
\begin{eqnarray*}
Q_n=\frac{v_n}{v_{n-1}}\left(1-\gamma_n\right)^mQ_{n-1}+\gamma_n\alpha_n
\end{eqnarray*}
with, since $\left(v_n\right)\in \mathcal{GS}\left(v^*\right)$ and in view of~\eqref{eq:6},
\begin{eqnarray}\label{eq:16}
\frac{v_n}{v_{n-1}}\left(1-\gamma_n\right)^m&=&\left(1+\frac{v^*}{n}+o\left(\frac{1}{n}\right)\right)\left(1-m\gamma_n+o\left(\gamma_n\right)\right)\nonumber\\
&=&\left(1+v^*\xi\gamma_n+o\left(\gamma_n\right)\right)\left(1-m\gamma_n+o\left(\gamma_n\right)\right)\nonumber\\
&=&1-\left(m-v^*\xi\right)\gamma_n+o\left(\gamma_n\right).
\end{eqnarray}
Set $A\in\left]0,m-v^*\xi\right[$; for $n$ large enough, we obtain 
\begin{eqnarray*}
Q_n&\leq&\left(1-A\gamma_n\right)Q_{n-1}+\gamma_n\alpha_n
\end{eqnarray*}
and~\eqref{eq:15} follows straightforwardly from the application of Lemma 4.I.1 in Duflo (1996). Now, let $C$ denote a positive generic constant that may vary from line to line; we have
\begin{eqnarray*}
v_n\Pi_n^{m}\sum_{k=1}^n\Pi_k^{-m}\gamma_kv_k^{-1}-\left(m-v^*\xi\right)^{-1}   &=&v_n\Pi_n^{m}\left[\sum_{k=1}^n\Pi_k^{-m}\gamma_kv_k^{-1}-\left(m-v^*\xi\right)^{-1}P_n\right]
\end{eqnarray*}
with, in view of~\eqref{eq:16},
\begin{eqnarray*}
P_n&=& v_n^{-1}\Pi_n^{-m}\\
&=&\sum_{k=2}^n\left(v_k^{-1}\Pi_k^{-m}-v_{k-1}^{-1}\Pi_{k-1}^{-m}\right)+C\\
&=&\sum_{k=2}^nv_k^{-1}\Pi_k^{-m}\left[1-\frac{v_k}{v_{k-1}}\left(1-\gamma_k\right)^m\right]+C\\
&=&\sum_{k=2}^nv_k^{-1}\Pi_k^{-m}\left[\left(m-v^*\xi\right)\gamma_k+o\left(\gamma_k\right)\right]+C.
\end{eqnarray*}
It follows that
\begin{eqnarray*}
v_n\Pi_n^{m}\sum_{k=1}^n\Pi_k^{-m}\gamma_kv_k^{-1}-\left(m-v^*\xi\right)^{-1}=v_n\Pi_n^{m}\left[\sum_{k=1}^n\Pi_k^{-m}v_k^{-1}o\left(\gamma_k\right)+C\right],
\end{eqnarray*}
and~\eqref{eq:14} follows from the application of~\eqref{eq:15}, which 
concludes the proof of Lemma~\ref{L:1}.

\subsection{Proof of Proposition~\ref{Pr:1}} \label{Section 5.2}
In view of~\eqref{eq:1} and~\eqref{eq:13}, we have
\begin{eqnarray}
\lefteqn{f_n\left(x\right)-f\left(x\right)}\nonumber\\
&=&\left(1-\gamma_n\right)\left(f_{n-1}\left(x\right)-f\left(x\right)\right)+\gamma_n\left(Z_n\left(x\right)-f\left(x\right)\right)\nonumber\\
&=&\sum_{k=1}^{n-1}\left[\prod_{j=k+1}^{n}\left(1-\gamma_j\right)\right]\gamma_k\left(Z_k\left(x\right)-f\left(x\right)\right)+\gamma_n\left(Z_n\left(x\right)-f\left(x\right)\right)+\left[\prod_{j=1}^{n}\left(1-\gamma_j\right)\right]\left(f_{0}\left(x\right)-f\left(x\right)\right)\nonumber\\
&=&\Pi_n\sum_{k=1}^n\Pi_k^{-1}\gamma_k\left(Z_k\left(x\right)-f\left(x\right)\right)+\Pi_n\left(f_0\left(x\right)-f\left(x\right)\right).
\label{encore}
\end{eqnarray}
It follows that
\begin{eqnarray*}
\mathbb{E}\left(f_n\left(x\right)\right)-f\left(x\right)&=&\Pi_n\sum_{k=1}^{n}\Pi_k^{-1}\gamma_k\left(\mathbb{E}\left(Z_k\left(x\right)\right)-f\left(x\right)\right)+\Pi_n\left(f_0\left(x\right)-f\left(x\right)\right).
\end{eqnarray*}
Taylor's expansion with integral remainder ensures that
\begin{eqnarray}
\mathbb{E}\left[Z_k\left(x\right)\right]-f\left(x\right)
&=&\int_{\mathbb{R}^d}K\left(z\right)\left[f\left(x-zh_k\right)-f\left(x\right)\right]dz
\nonumber\\
&=&\frac{h_k^2}{2}
\sum_{j=1}^d\left(\mu_j^2f^{(2)}_{jj}(x)\right)
+h_k^2\delta_k\left(x\right)
\label{eq:17}
\end{eqnarray}
with
\begin{eqnarray*}
\delta_k\left(x\right)=\sum_{1\leq i,j\leq d}\int_{\mathbb{R}^d}
\int_{0}^{1}\left(1-s\right)z_iz_jK(z)\left[f^{(2)}_{ij}\left(x-zh_ks\right)
-f^{\left(2\right)}_{ij}\left(x\right)\right]dsdz,
\end{eqnarray*}
and, since $f^{\left(2\right)}_{ij}$ is bounded and continuous at $x$ for 
all $i,j\in\{1,\ldots ,d\}$, we have 
$\lim_{k\to\infty}\delta_k\left(x\right)=0$. 
In the case $a\leq\al/(d+4)$, we have $\lim_{n\to\infty}\left(n\gamma_n\right)>2a$; 
the application of Lemma $\ref{L:1}$ then gives 
\begin{eqnarray*}
\mathbb{E}\left[f_n\left(x\right)\right]-f\left(x\right)&=&
\frac{1}{2}
\sum_{j=1}^d\left(\mu_j^2f^{(2)}_{jj}(x)\right)
\Pi_n\sum_{k=1}^{n}\Pi_k^{-1}\gamma_kh_k^2[1+o(1)]+\Pi_n\left(f_0\left(x\right)-f\left(x\right)\right) \nonumber\\
&=&
\frac{1}{2(1-2a\xi)}
\sum_{j=1}^d\left(\mu_j^2f^{(2)}_{jj}(x)\right)
\left[h_n^2+o(1)\right],
\end{eqnarray*}
and~\eqref{eq:7} follows.
In the case $a>\al/(d+4)$, we have $h_n^2=o\left(\sqrt{\gamma_nh_n^{-d}}\right)$; 
since $\lim_{n\to\infty}\left(n\gamma_n\right)>\left(\alpha-ad\right)/2$, 
Lemma $\ref{L:1}$ then ensures that 
\begin{eqnarray*}
\mathbb{E}\left[f_n\left(x\right)\right]-f\left(x\right)
&=&\Pi_n\sum_{k=1}^{n}\Pi_k^{-1}\gamma_ko\left(\sqrt{\gamma_kh_k^{-d}}\right)
+O\left(\Pi_n\right)\\
&=&o\left(\sqrt{\gamma_nh_n^{-d}}\right),
\end{eqnarray*}
which gives \eqref{eq:8}. Now, we have
\begin{eqnarray*}
Var\left[f_n\left(x\right)\right]
&=&\Pi_n^{2}\sum_{k=1}^n\Pi_k^{-2}\gamma_k^2Var\left[Z_k\left(x\right)\right]\nonumber\\
&=& \Pi_n^{2}\sum_{k=1}^n\frac{\Pi_k^{-2}\gamma_k^2}{h_k^d}
\left[\int_{\mathbb{R}^d}K^2\left(z\right)f\left(x-zh_k\right)dz-
h_k^d\left(\int_{\mathbb{R}^d}K\left(z\right)f\left(x-zh_k\right)dz\right)^2\right]\\
&=& \Pi_n^{2}\sum_{k=1}^n\frac{\Pi_k^{-2}\gamma_k^2}{h_k^d}
\left[f\left(x\right)\int_{\mathbb{R}^d}K^2\left(z\right)dz+
\nu_k\left(x\right)-h_k^d\tilde \nu_k\left(x\right)\right]
\end{eqnarray*}
with
\begin{eqnarray*}
\nu_k\left(x\right)&=&
\int_{\mathbb{R}^d}K^2\left(z\right)\left[f\left(x-zh_k\right)-f\left(x\right)\right]dz,\\
\tilde \nu_k\left(x\right)&=&
\left(\int_{\mathbb{R}^d}K\left(z\right)f\left(x-zh_k\right)dz\right)^2.
\end{eqnarray*}
Since $f$ is bounded and continuous, we have $\lim_{k\to\infty}\nu_k\left(x\right)=0$ and 
$\lim_{k\to\infty}h_k^d\tilde\nu_k\left(x\right)=0$. 
In the case $a\geq \al/(d+4)$, we have 
$\lim_{n\to\infty}\left(n\gamma_n\right)>\left(\alpha-ad\right)/2 $, 
and the application of Lemma~\ref{L:1} gives
\begin{eqnarray*}
Var\left[f_n\left(x\right)\right]
&=& \Pi_n^{2}\sum_{k=1}^n\frac{\Pi_k^{-2}\gamma_k^2}{h_k^d}\left[
f\left(x\right)\int_{\mathbb{R}^d}K^2\left(z\right)dz+o\left(1\right)\right]\\
&=& \frac{1}{2-\left(\alpha-ad\right)\xi}\frac{\gamma_n}{h_n^d}
\left[
f\left(x\right)\int_{\mathbb{R}^d}K^2\left(z\right)dz+o\left(1\right)\right],
\end{eqnarray*}
which proves~\eqref{eq:9}. 
In the case $a<\al/(d+4)$, we have $\gamma_nh_n^{-d}=o\left(h_n^4\right)$; since 
$\lim_{n\to\infty}\left(n\gamma_n\right)>2a$, 
Lemma~\ref{L:1} then ensures that 
\begin{eqnarray*}
Var\left[f_n\left(x\right)\right]
&=& \Pi_n^{2}\sum_{k=1}^n\Pi_k^{-2}\gamma_ko\left(h_k^4\right)\\
&=& o\left(h_n^4\right),
\end{eqnarray*}
which gives~\eqref{eq:10}.

\subsection{Proof of Proposition~\ref{Pr:2}} \label{Section 5.3}
Let us first note that, in view of \eqref{eq:17}, we have 
\begin{eqnarray*}
\lefteqn{\int_{\mathbb{R}^d}\left\{
\Pi_n\sum_{k=1}^{n}\Pi_k^{-1}\gamma_k\left[\mathbb{E}\left(Z_k\left(x\right)\right)-f\left(x\right)\right]
\right\}^2dx}\\
&=&\frac{1}{4}\int_{\mathbb{R}^d}\left[\sum_{j=1}^d\mu_j^2f^{(2)}_{jj}(x)\right]^2
dx\left[\Pi_n\sum_{k=1}^{n}\Pi_k^{-1}\gamma_kh_k^2\right]^2
+\int_{\mathbb{R}^d}\left[\Pi_n\sum_{k=1}^n\Pi_k^{-1}\gamma_kh_k^2\delta_k\left(x\right)\right]^2dx
\nonumber\\&&
+\left(\Pi_n\sum_{k=1}^{n}\Pi_k^{-1}\gamma_kh_k^2\right)\left(\Pi_n\sum_{k=1}^{n}\Pi_k^{-1}\gamma_kh_k^2\int_{\mathbb{R}^d}\left[\sum_{j=1}^d\mu_j^2f^{(2)}_{jj}(x)\right]
\delta_k\left(x\right)dx\right).
\end{eqnarray*}
Since $f^{\left(2\right)}_{ij}$ is continuous, bounded, and integrable for 
all $i,j\in\{1,\ldots d\}$, the application 
of Lebesgue's convergence theorem ensures that $\lim_{k \to +\infty}\int_{\mathbb{R}^d}\delta_k^2\left(x\right)dx=0$ and 
$\lim_{k \to +\infty}\int_{\mathbb{R}^d}
[\sum_{j=1}^d\mu_j^2f^{(2)}_{jj}(x)]\delta_k\left(x\right)dx=0$.
Moreover, Jensen's inequality gives
\begin{eqnarray*}
\int_{\mathbb{R}^d}\left[\Pi_n\sum_{k=1}^n\Pi_k^{-1}\gamma_kh_k^2\delta_k\left(x\right)\right]^2dx
&\leq &\left(\Pi_n\sum_{k=1}^n\Pi_k^{-1}\gamma_kh_k^2\right)\left(\Pi_n\sum_{k=1}^{n}\Pi_k^{-1}\gamma_kh_k^2\int_{\mathbb{R}^d}\delta_k^2\left(x\right)dx\right)\nonumber\\
& \leq & \left(\Pi_n\sum_{k=1}^n\Pi_k^{-1}\gamma_kh_k^2\right)\left(\Pi_n\sum_{k=1}^{n}\Pi_k^{-1}\gamma_ko\left(h_k^2\right)\right),\nonumber\\
\end{eqnarray*}
so that we get 
\begin{eqnarray*}
\lefteqn{\int_{\mathbb{R}^d}\left\{
\Pi_n\sum_{k=1}^{n}\Pi_k^{-1}\gamma_k\left[\mathbb{E}\left(Z_k\left(x\right)\right)-f\left(x\right)\right]
\right\}^2dx}\\
&=&\frac{1}{4}\int_{\mathbb{R}^d}\left[\sum_{j=1}^d\mu_j^2f^{(2)}_{jj}(x)\right]^2
dx\left[\Pi_n\sum_{k=1}^{n}\Pi_k^{-1}\gamma_kh_k^2\right]^2
+O\left(\left[\Pi_n\sum_{k=1}^{n}\Pi_k^{-1}\gamma_kh_k^2\right]
\left[\Pi_n\sum_{k=1}^{n}\Pi_k^{-1}\gamma_ko(h_k^2)
\right]\right).
\end{eqnarray*}
$\bullet$ Let us first consider the case $a\leq\alpha/(d+4)$. In this case, 
$\lim_{n\to\infty}\left(n\gamma_n\right)>2a$, and 
the application of Lemma $\ref{L:1}$ gives 
\begin{eqnarray*}
\int_{\mathbb{R}^d}\left\{
\Pi_n\sum_{k=1}^{n}\Pi_k^{-1}\gamma_k\left[\mathbb{E}\left(Z_k\left(x\right)\right)-f\left(x\right)\right]
\right\}^2dx
& = & \frac{1}{4\left(1-2a\xi\right)^2}h_n^4\int_{\mathbb{R}^d}
\left[\sum_{j=1}^d\mu_j^2f^{(2)}_{jj}(x)\right]^2dx+o\left(h_n^4\right),
\end{eqnarray*}
and ensures that $\Pi_n^2=o(h_n^4)$. In view of (\ref{encore}), we then deduce that
\begin{eqnarray}\label{eq:609}
\int_{\mathbb{R}^d}\left\{\mathbb{E}\left( f_{n}\left( x\right) \right) -f\left( x\right)\right\}^2dx=\frac{1}{4\left(1-2a\xi\right)^2}h_n^4\int_{\mathbb{R}^d}
\left[\sum_{j=1}^d\mu_j^2f^{(2)}_{jj}(x)\right]^2dx+o\left(h_n^4\right).
\end{eqnarray}
$\bullet$ Let us now consider the case $a>\alpha/(d+4)$. 
In this case, we have $h_k^2=o(\sqrt{\g_kh_k^{-d}})$ and 
$\lim_{n\to\infty}\left(n\gamma_n\right)>\left(\alpha-ad\right)/2$.   
The application of Lemma $\ref{L:1}$ then gives 
\begin{eqnarray*}
\int_{\mathbb{R}^d}\left\{
\Pi_n\sum_{k=1}^{n}\Pi_k^{-1}\gamma_k\left[\mathbb{E}\left(Z_k\left(x\right)\right)-f\left(x\right)\right]
\right\}^2dx
&=&O\left(\left[\Pi_n\sum_{k=1}^{n}\Pi_k^{-1}\gamma_ko\left(\sqrt{\g_kh_k^{-d}}\right)
\right]^2\right)\\
& = & 
o\left(\g_nh_n^{-d}\right),
\end{eqnarray*}
and ensures that $\Pi_n^2=o(\g_nh_n^{-d})$. In view of (\ref{encore}), we then deduce that
\begin{eqnarray}\label{biseq:609}
\int_{\mathbb{R}^d}\left\{\mathbb{E}\left( f_{n}\left( x\right) \right) -f\left( x\right)\right\}^2dx=
o\left(\g_nh_n^{-d}\right).
\end{eqnarray}
On the other hand, we note that 
\begin{eqnarray*}
\lefteqn{\int_{\mathbb{R}^d}Var\left[f_n\left(x\right)\right]dx}\nonumber\\
&=&\Pi_n^{2}\sum_{k=1}^n\Pi_k^{-2}\gamma_k^2\int_{\mathbb{R}^d}Var\left[Z_k\left(x\right)\right]dx\nonumber\\
&=& \Pi_n^{2}\sum_{k=1}^n\Pi_k^{-2}\gamma_k^2
\left[\frac{1}{h_k^d}\int_{\mathbb{R}^d}\int_{\mathbb{R}^d}K^2\left(z\right)
f\left(x-zh_k\right)dzdx- 
\int_{\mathbb{R}^d}\left(\int_{\mathbb{R}^d}K\left(z\right)f\left(x-zh_k\right)dz\right)^2dx
\right]
\end{eqnarray*}
with 
\begin{eqnarray*}
\int_{\mathbb{R}^d}\int_{\mathbb{R}^d}K^2\left(z\right)f\left(x-zh_k\right)dzdx&=&\int_{\mathbb{R}^d}K^2\left(z\right)\left(\int_{\mathbb{R}^d}f\left(x-zh_k\right)dx\right)dz\\
&=&\int_{\mathbb{R}^d}K^2\left(z\right)dz
\end{eqnarray*}
and
\begin{eqnarray*}
\int_{\mathbb{R}^d}\left(\int_{\mathbb{R}^d}K\left(z\right)f\left(x-zh_k\right)dz\right)^2dx
&=&\int_{\mathbb{R}^{3d}}K\left(z\right)K\left(z'\right)f\left(x-zh_k\right)f\left(x-z'h_k\right)dz dz'dx\\
&\leq &\|f\|_{\infty}\|K\|_1^2.
\end{eqnarray*}
$\bullet$
In the case $a\geq\alpha/(d+4)$, we have 
$\lim_{n\to\infty}\left(n\gamma_n\right)>\left(\alpha-ad\right)/2$, 
and Lemma $\ref{L:1}$ ensures that
\begin{eqnarray}\label{eq:6/09}
\int_{\mathbb{R}^d}Var\left[f_n\left(x\right)\right]dx
&=&\Pi_n^{2}\sum_{k=1}^n\frac{\Pi_k^{-2}\gamma_k^2}{h_k^d}\left[
\int_{\mathbb{R}^d}K^2\left(z\right)dz+o(1)\right]\nonumber\\
&=&\frac{\gamma_n}{h_n^d}\frac{1}{(2-\left(\alpha-ad\right)\xi)}
\int_{\mathbb{R}^d}K^2\left(z\right)dz+o\left(\frac{\gamma_n}{h_n^d}\right).
\end{eqnarray}
$\bullet$
In the case $a<\alpha/(d+4)$, we have $\gamma_nh_n^{-d}=o(h_n^4)$ and 
$\lim_{n\to\infty}\left(n\gamma_n\right)>2a$, so that 
Lemma $\ref{L:1}$ gives 
\begin{eqnarray}\label{biseq:6/09}
\int_{\mathbb{R}^d}Var\left[f_n\left(x\right)\right]dx
&=&\Pi_n^{2}\sum_{k=1}^n\Pi_k^{-2}\gamma_ko\left(h_k^4\right) \nonumber\\
&=&o\left(h_n^4\right).
\end{eqnarray}
Part 1 of Proposition~\ref{Pr:2} follows from the combination of~\eqref{eq:609} 
and \eqref{biseq:6/09}, Part 2 from that of~\eqref{eq:609} and~\eqref{eq:6/09}, 
and Part 3 from that of~\eqref{biseq:609} and~\eqref{eq:6/09}.

\subsection{Proof of Theorem~\ref{T:1}} \label{Section 5.4}
Let us at first assume that, if $a\geq\alpha/(d+4)$, then 
\begin{eqnarray}\label{eq:22}
\sqrt{\gamma_n^{-1} h_n^d}\left(f_{n}\left( x\right)-\mathbb{E}\left[f_n\left(x\right)\right]\right) \stackrel{\mathcal{D}}{\rightarrow}\mathcal{N}\left( 0,
\frac{1}{2-\left(\alpha-ad\right)\xi}f\left(x\right) 
\int_{\mathbb{R}^d}K^2\left(z\right)dz\right).
\end{eqnarray}
In the case when $a>\alpha/(d+4)$, Part 1 of Theorem~\ref{T:1} follows from the combination of~\eqref{eq:8} and~\eqref{eq:22}. In the case when $a=\alpha/(d+4)$, Parts 1 and 2 of Theorem~\ref{T:1} follow from the combination of~\eqref{eq:7} and~\eqref{eq:22}. In the case $a<\alpha/(d+4)$,~\eqref{eq:10} implies that 
\begin{eqnarray*}
h_n^{-2}\left(f_n\left(x\right)-\mathbb{E}\left(f_n\left(x\right)\right)\right)\stackrel{\mathbb{P}}{\rightarrow}0,
\end{eqnarray*}
and the application of~\eqref{eq:7} gives Part 2 of Theorem~\ref{T:1}.\\

We now prove~\eqref{eq:22}. 
In view of~\eqref{eq:1}, we have
\begin{eqnarray*}
f_n\left(x\right)-\mathbb{E}\left[f_n\left(x\right)\right]
&=&\left(1-\gamma_n\right)\left(f_{n-1}\left(x\right)-
\mathbb{E}\left[f_{n-1}\left(x\right)\right]
\right)+\gamma_n\left(Z_n\left(x\right)-\mathbb{E}\left[Z_n\left(x\right)\right]
\right)\\
&=&\Pi_n\sum_{k=1}^n\Pi_k^{-1}\gamma_k\left(Z_k\left(x\right)-
\mathbb{E}\left[Z_k\left(x\right)\right]\right).
\end{eqnarray*}
Set 
\begin{eqnarray}\label{eq:23}
Y_{k}\left(x\right)=\Pi_k^{-1}\gamma_k\left(Z_k\left(x\right)-\mathbb{E}\left(Z_k\left(x\right)\right)\right).
\end{eqnarray}
The application of Lemma $\ref{L:1}$ ensures that
\begin{eqnarray}\label{eq:24} 
v_n^2&=&\sum_{k=1}^nVar\left(Y_{k}\left(x\right)\right)\nonumber\\
&=&\sum_{k=1}^n\Pi_k^{-2}\gamma_k^2Var\left(Z_k\left(x\right)\right)\nonumber\\
&=&\sum_{k=1}^n\frac{\Pi_k^{-2}\gamma_k^2}{h_k^d}\left[f\left(x\right)\int_{\mathbb{R}^d}K^2\left(z\right)dz+o\left(1\right)\right]\nonumber\\
&=&\frac{1}{\Pi_n^2}\frac{\gamma_n}{h_n^d}\left[\frac{1}{2-\left(\alpha-ad\right)\xi}f\left(x\right)\int_{\mathbb{R}^d}K^2\left(z\right)dz+o\left(1\right)\right].
\end{eqnarray}
On the other hand, we have, for all $p>0$, 
\begin{eqnarray}
\label{eq:25}
\mathbb{E}\left[\left|Z_k\left(x\right)\right|^{2+p}\right] &=&
O\left(\frac{1}{h_k^{d(1+p)}}\right),
\end{eqnarray}
and, since $\lim_{n\to\infty}\left(n\gamma_n\right)>\left(\alpha-ad\right)/2$, there exists $p>0$ such that $\lim_{n\to \infty}\left(n\gamma_n\right)>\frac{1+p}{2+p}\left(\alpha-ad\right)$. Applying Lemma $\ref{L:1}$, we get 
\begin{eqnarray*}
\sum_{k=1}^n\mathbb{E}\left[\left|Y_{k}\left(x\right)\right|^{2+p}\right]&=&O\left(\sum_{k=1}^n \Pi_k^{-2-p}\gamma_k^{2+p}\mathbb{E}\left[\left|Z_k\left(x\right)\right|^{2+p}\right]\right)\nonumber\\
&=&O\left(\sum_{k=1}^n \frac{\Pi_k^{-2-p}\gamma_k^{2+p}}{h_k^{d(1+p)}}\right)\\
&=&O\left(\frac{\gamma_n^{1+p}}{\Pi_n^{2+p}h_n^{d(1+p)}}\right)\nonumber,
\end{eqnarray*}
and we thus obtain 
\begin{eqnarray*}
\frac{1}{v_n^{2+p}}\sum_{k=1}^n\mathbb{E}\left[\left|Y_{k}\left(x\right)\right|^{2+p}\right]& = & O\left({\left[\gamma_nh_n^{-d}\right]}^{p/2}\right)=o\left(1\right).
\end{eqnarray*}
The convergence in~\eqref{eq:22} then follows from the application of Lyapounov's Theorem.

\subsection{Proof of Theorem~\ref{T:2}} \label{Section 5.5}
Set
\begin{eqnarray*}
S_n\left(x\right)&=&\sum_{k=1}^nY_{k}\left(x\right)
\end{eqnarray*}
where $Y_k$ is defined in~\eqref{eq:23}, and set $\gamma_0=h_0=1$.\\
$\bullet$ Let us first consider the case $a\geq \alpha/(d+4)$ (in which case $\lim_{n\to\infty}\left(n\gamma_n\right)>\left(\alpha-ad\right)/2$). We set $H_n^2=\Pi_n^2\gamma_n^{-1}h_n^d$, and note that, since $\left(\gamma_n^{-1}h_n^d\right)\in \mathcal{GS}\left(\alpha-ad\right)$, we have
\begin{eqnarray}\label{eq:26}
\ln \left(H_n^{-2}\right)&=&-2\ln \left(\Pi_n\right)+\ln\left(\prod_{k=1}^n\frac{\gamma_{k-1}^{-1}h_{k-1}^d}{\gamma_k^{-1}h_k^d}\right)\nonumber\\
&=&-2\sum_{k=1}^n\ln\left(1-\gamma_k\right)+\sum_{k=1}^n\ln\left(1-\frac{\alpha-ad}{k}+o\left(\frac{1}{k}\right)\right)\nonumber\\
&=&\sum_{k=1}^n\left(2\gamma_k+o\left(\gamma_k\right)\right)-\sum_{k=1}^n\left(\left(\alpha-ad\right)\xi\gamma_k+o\left(\gamma_k\right)\right)\nonumber\\
&=&\left(2-\xi\left(\alpha-ad\right)\right)s_n+o\left(s_n\right).
\end{eqnarray}
Since $2-\xi\left(\alpha-ad\right)>0$, it follows in particular that $\lim_{n\to+\infty}H_n^{-2}=\infty$. Moreover, we clearly have $\lim_{n\to+\infty}H_n^2/H_{n-1}^2=1$, and by~\eqref{eq:24} 
\begin{eqnarray*}
\lim_{n\to+\infty}H_n^2\sum_{k=1}^nVar\left[Y_k\left(x\right)\right]=\frac{1}
{2-\left(\alpha-ad\right)\xi}f\left(x\right)\int_{\mathbb{R}^d}K^2\left(z\right)dz.
\end{eqnarray*}
Now, in view of~\eqref{eq:25}, $\mathbb{E}\left[\left|Y_{k}\left(x\right)\right|^{3}\right]=O\left(\Pi_k^{-3}\gamma_k^3h_k^{-2d}\right)$ and, since $\lim_{n\to\infty}\left(n\gamma_n\right)>\left(\alpha-ad\right)/2$, the application of Lemma~\ref{L:1} and of~\eqref{eq:26} gives
\begin{eqnarray*}
\frac{1}{n\sqrt{n}}\sum_{k=1}^n\mathbb{E}\left(\left|H_nY_k\left(x\right)\right|^3\right)&=&O\left(\frac{H_n^3}{n\sqrt{n}}\sum_{k=1}^n \Pi_k^{-3}\gamma_k^{3}h_k^{-2d}\right)\\
&=&O\left(\frac{H_n^3}{n\sqrt{n}}\sum_{k=1}^n \Pi_k^{-3}\gamma_k
o\left(\left[\gamma_kh_k^{-d}\right]^{{3}/{2}}\right)\right)\\
&=&o\left(\frac{H_n^3}{n\sqrt{n}}\Pi_n^{-3}\left[\gamma_nh_n^{-d}\right]^{{3}/{2}}
\right)\\
&=&o\left(\frac{1}{n\sqrt{n}}\right)\\
&=&o\left(\left[\ln \left(H_n^{-2}\right)\right]^{-1}\right).
\end{eqnarray*}
The application of Theorem 1 in Mokkadem and Pelletier (2007b) then ensures that, with probability one, the sequence
\begin{eqnarray*}
\left(\frac{H_nS_n\left(x\right)}{\sqrt{2\ln \ln \left(H_n^{-2}\right)}}\right)=\left(\frac{\sqrt{\gamma_n^{-1}h_n^d}\left(f_{n}\left( x\right)-\mathbb{E}\left(f_n\left(x\right)\right)\right)}{\sqrt{2\ln \ln \left(H_n^{-2}\right)}}\right)
\end{eqnarray*}
is relatively compact and its limit set is the interval 
\begin{eqnarray}\label{eq:27}
\left[-\sqrt{\frac{f\left(x\right)}{2-\left(\alpha-ad\right)\xi}
\int_{\mathbb{R}^d} K^2\left(z\right)dz},\sqrt{\frac{f\left(x\right)}
{2-\left(\alpha-ad\right)\xi}
\int_{\mathbb{R}^d} K^2\left(z\right)dz}\right].
\end{eqnarray}
In view of~\eqref{eq:26}, we have $\lim_{n\to\infty}\ln \ln \left(H_n^{-2}\right)/\ln s_n=1$. It follows that, with probability one, the sequence $\left(\sqrt{\gamma_n^{-1}h_n^d}\left(f_{n}\left( x\right)-\mathbb{E}\left(f_n\left(x\right)\right)\right)/\sqrt{2\ln s_n}\right)$ is relatively compact, and its limit set is the interval given in~\eqref{eq:27}. The application of~\eqref{eq:7} (respectively~\eqref{eq:8}) concludes the proof of Theorem~\ref{T:2} in the case $a=\alpha/(d+4)$ (respectively $a>\alpha/(d+4)$).\\
$\bullet$ Let us now consider the case $a<\alpha/(d+4)$ (in which case $\lim_{n\to\infty}\left(n\gamma_n\right)>2a$). Set $H_n^{-2}=\Pi_n^{-2}h_n^{4}\left(\ln \ln \left(\Pi_n^{-2}h_n^4\right)\right)^{-1}$, and note that, since $\left(h_n^{-4}\right)\in \mathcal{GS}\left(4a\right)$, we have
\begin{eqnarray}\label{eq:28}
\ln \left(\Pi_n^{-2}h_n^{4}\right)&=&-2\ln \left(\Pi_n\right)+\ln\left(\prod_{k=1}^n\frac{h_{k-1}^{-4}}{h_k^{-4}}\right)\nonumber\\
&=&-2\sum_{k=1}^n\ln\left(1-\gamma_k\right)+\sum_{k=1}^n\ln\left(1-\frac{4a}{k}+o\left(\frac{1}{k}\right)\right)\nonumber\\
&=&\sum_{k=1}^n\left(2\gamma_k+o\left(\gamma_k\right)\right)-\sum_{k=1}^n\left(4a\xi\gamma_k+o\left(\gamma_k\right)\right)\nonumber\\
&=&\left(2-4a\xi\right)s_n+o\left(s_n\right).
\end{eqnarray}
Since $2-4a\xi>0$, it follows in particular that $\lim_{n\to\infty}\Pi_n^{-2}h_n^4=\infty$, and thus $\lim_{n\to\infty}H_n^{-2}=\infty$. Moreover, we clearly have $\lim_{n\to\infty}H_n^2/H_{n-1}^2=1$. Set $\epsilon \in \left]0,\alpha-(d+4)a\right[$ such that $\lim_{n\to\infty}\left(n\gamma_n\right)>2a+\epsilon/2$; in view of~\eqref{eq:24}, and applying Lemma~\ref{L:1}, we get 
\begin{eqnarray*}
H_n^2\sum_{k=1}^nVar\left[Y_{k}\left(x\right)\right]&=&O\left(\Pi_n^2h_n^{-4}\ln \ln \left(\Pi_n^{-2}h_n^{4}\right)\sum_{k=1}^n\frac{\Pi_k^{-2}\gamma_k^2}{h_k^d}\right)\\
&=&O\left(\Pi_n^2h_n^{-4}\ln \ln \left(\Pi_n^{-2}h_n^{4}\right)\sum_{k=1}^n\Pi_k^{-2}\gamma_ko\left(h_k^4k^{-\epsilon}\right)\right)\\
&=&o\left(\ln \ln\left(\Pi_n^{-2}h_n^4\right)n^{-\epsilon}\right)\\
&=&o\left(1\right).
\end{eqnarray*}
Moreover, applying~\eqref{eq:25}, Lemma~\ref{L:1}, and~\eqref{eq:28}, we obtain
\begin{eqnarray*}
\frac{1}{n\sqrt{n}}\sum_{k=1}^n\mathbb{E}\left(\left|H_nY_k\left(x\right)\right|^3\right)&=&O\left(\frac{\Pi_n^3h_n^{-6}}{n\sqrt{n}}\left[\ln \ln \left(\Pi_n^{-2}h_n^4\right)\right]^{\frac{3}{2}}\left(\sum_{k=1}^n \Pi_k^{-3}\gamma_k^{3}h_k^{-2d}\right)\right)\\
&=&O\left(\frac{\Pi_n^3h_n^{-6}}{n\sqrt{n}}\left[\ln \ln \left(\Pi_n^{-2}h_n^4\right)\right]^{\frac{3}{2}}\left(\sum_{k=1}^n \Pi_k^{-3}\gamma_ko\left(h_k^{6}\right)\right)\right)\\
&=&o\left(\frac{\Pi_n^3h_n^{-6}}{n\sqrt{n}}\Pi_n^{-3}h_n^{6}\left[\ln \ln \left(\Pi_n^{-2}h_n^4\right)\right]^{\frac{3}{2}}\right)\\
&=&o\left(\left[\ln \left(H_n^{-2}\right)\right]^{-1}\right).
\end{eqnarray*}
The application of Theorem 1 in Mokkadem and Pelletier (2007b) then ensures that, with probability one, 
\begin{eqnarray*}
\lim_{n\to\infty}\frac{H_nS_n\left(x\right)}{\sqrt{2\ln \ln \left(H_n^{-2}\right)}}=\lim_{n\to\infty}h_n^{-2}\frac{\sqrt{\ln \ln \left(\Pi_n^{-2}h_n^4\right)}}{\sqrt{2\ln \ln \left(H_n^{-2}\right)}}\left(f_n\left(x\right)-\mathbb{E}\left(f_n\left(x\right)\right)\right)=0.
\end{eqnarray*}
Noting that~\eqref{eq:28} ensures that $\lim_{n\to\infty}\ln \ln \left(H_n^{-2}\right)/\ln \ln \left(\Pi_n^{-2}h_n^4\right)=1$, we deduce that
\begin{eqnarray*}
\lim_{n\to\infty}h_n^{-2}\left[T_{n}\left( x\right)-\mathbb{E}\left(T_n\left(x\right)\right)\right]=0\quad a.s.,
\end{eqnarray*}
and Theorem~\ref{T:2} in the case $a<\alpha/(d+4)$ follows from~\eqref{eq:7}.

\subsection{Proof of Corollaries \ref{corollaire 3}-\ref{corollaire niveau}} \label{Section 5.2bis}

In view of~\eqref{eq:9}, 
to minimize the variance of $f_n$, the stepsize $\left(\gamma_n\right)$ must belong to $\mathcal{GS}\left(-1\right)$ and satisfy $\lim_{n\to\infty}n\gamma_n=\gamma_0\in]0,\infty[$. For such a choice, $\xi=\gamma_0^{-1}$, so that~\eqref{eq:9} can be rewritten as
$$
Var\left( f_{n}\left( x\right) \right) =\frac{\gamma_0}{2-\left(1-ad\right)\gamma_0^{-1}}\frac{1}{nh_n^d}f\left(x\right)\int_{\mathbb{R}^d} K^2\left(z\right)dz+o\left(\frac{1}{nh_n^d}\right). 
$$
The function $\gamma_0 \mapsto \gamma_0\left[2-\left(1-ad\right)\gamma_0^{-1}\right]^{-1} $ reaching its minimum at the point $\gamma_0=1-ad$, Corollary \ref{corollaire 3} follows.\\

Let us now prove Corollary \ref{corollaire niveau}. When $\lim_{n\to \infty}n\gamma_n=\g_0>0$ 
and $\lim_{n\to \infty}nh_n^d=0$, the first part of Theorem \ref{T:1} ensures that 
\begin{eqnarray*}
\sqrt{nh_{n}^d}\left(f_{n}\left( x\right)-f\left( x\right) \right) 
& \stackrel{\mathcal{D}}{\rightarrow} & \mathcal{N}\left(0,\frac{\g_0^2}{2\g_0-(1-ad)}f\left(x\right) 
\int_{\mathbb{R}^d}K^2\left(z\right)dz\right).
\end{eqnarray*}
Proposition \ref{Pr:1} ensuring the consistency of $f_n$, Corollary \ref{corollaire niveau} follows. 
\\

We now show how Corollary~\ref{C:1} can be deduced from Proposition~\ref{Pr:1}. Corollary~\ref{C:2} 
is deduced from Proposition~\ref{Pr:2} exactly in the same way, so that its proof is omitted. Set
\begin{eqnarray*}
C_1\left(\xi\right)&=&\frac{1}{4\left(1-2a\xi\right)^2}
\left(\sum_{j=1}^d\mu_j^2f^{(2)}_{jj}(x)\right)^2,\\
C_2\left(\xi\right)&=&\frac{1}{2-\left(1-ad\right)\xi}f\left(x\right)\int_{\mathbb{R}^d} K^2\left(z\right)dz.
\end{eqnarray*}
The application of Proposition~\ref{Pr:1} ensures that
\begin{eqnarray}\label{eq:11}
MSE = \left\{\begin{array}{llll}C_1\left(\xi\right)h_n^4+o\left(h_n^{4}\right) & \mbox{ if } & a<\alpha/(d+4),\\
C_1\left(\xi\right)h_n^4+C_2\left(\xi\right)\gamma_nh_n^{-d}
+o\left(h_n^{4}+\gamma_nh_n^{-d}\right)  & \mbox{ if } & a=\alpha/(d+4),\\
C_2\left(\xi\right)\gamma_nh_n^{-d}+o\left(\gamma_nh_n^{-d}\right) & 
\mbox{ if } & a>\alpha/(d+4).
\end{array}\right.
\end{eqnarray}
Set $\alpha \in \left]1/2,1\right]$. If $a=\alpha/(d+4)$, $\left(C_1\left(\xi\right)h_n^4+C_2\left(\xi\right)\gamma_nh_n^{-d}\right)\in \mathcal{GS}\left(-4\alpha/(d+4)\right)$. If $a<\alpha/(d+4)$, $\left(h_n^4\right)\in \mathcal{GS}\left(-4a\right)$ with $-4a>-4\alpha/(d+4)$, and, if $a>\alpha/(d+4)$, $\left(\gamma_nh_n^{-d}\right)\in \mathcal{GS}\left(-\alpha+ad\right)$ with $-\alpha+ad>-4\alpha/(d+4)$. It follows that, for a given $\alpha$, to minimize the MSE of $f_n$, the parameter $a$ must be chosen equal to $\alpha/(d+4)$. Moreover, in view of~\eqref{eq:11}, the parameter $\alpha$ must be chosen equal to $1$. In other words, to minimize the MSE of $f_n$, the stepsize $\left(\gamma_n\right)$ must be chosen in $\mathcal{GS}\left(-1\right)$, the bandwidth $\left(h_n\right)$ in $\mathcal{GS}\left(-1/(d+4)\right)$ (and, in view of $(A2)iii)$, the condition $\lim_{n\to\infty}n\gamma_n>2/(d+4)$ must be fulfilled). For this choice of stepsize and bandwidth, set $\mathcal{L}_n=n\gamma_n$ and $\tilde{\mathcal{L}}_n=n^{1/(d+4)}h_n$. The MSE of $f_n$ can then be rewritten as
\begin{eqnarray*}
MSE=n^{-\frac{4}{d+4}}\left[C_1\left(\xi\right)\tilde{\mathcal{L}}_n^4+C_2\left(\xi\right)\mathcal{L}_n\tilde{\mathcal{L}}_n^{-d}\right]\left[1+o\left(1\right)\right].
\end{eqnarray*}
Now, set $\mathcal{L}_n$. Since the function 
$x \mapsto C_1\left(\xi\right)x^4+C_2\left(\xi\right)\mathcal{L}_nx^{-d}$ reaches 
its minimum at the point 
$\left(dC_2\left(\xi\right)\mathcal{L}_n/
\left[4C_1\left(\xi\right)\right]\right)^{1/(d+4)}$, 
to minimise the MSE of $f_n$, $\tilde{\mathcal{L}}_n$ must be chosen equal to 
$\left(dC_2\left(\xi\right)\mathcal{L}_n
/\left[4C_1\left(\xi\right)\right]\right)^{1/(d+4)}$, that is, 
$\left(h_n\right)$ must equal 
$\left(dC_2\left(\xi\right)/\left[4C_1\left(\xi\right)\right]\gamma_n\right)^{1/(d+4)}$. For such a choice, the MSE of $f_n$ can be rewritten as
\begin{eqnarray*}
MSE=n^{-\frac{4}{d+4}}\mathcal{L}_n^{\frac{4}{d+4}}
\left(\frac{d}{4}\right)^{-\frac{d}{d+4}}
\frac{d+4}{4}\left[C_1\left(\xi\right)\right]^{\frac{d}{d+4}}
\left[C_2\left(\xi\right)\right]^{\frac{4}{d+4}}\left[1+o\left(1\right)\right].
\end{eqnarray*}
It follows that to minimize the MSE of $f_n$, the limit of $\mathcal{L}_n$ 
(that is, of $n\gamma_n$) must be finite (and larger than $2/(d+4)$). 
Now, set $\gamma_0>2/(d+4)$ and $\mathcal{L}_n=\gamma_0\delta_n$ 
with $\lim_{n\to\infty}\delta_n=1$ (so that $\lim_{n\to\infty}n\gamma_n=\gamma_0$). 
In this case, we have $\xi=\gamma_0^{-1}$,
\begin{eqnarray*}
\begin{array}{llll}C_1\left(\xi\right)=
\frac{\gamma_0^2}{4\left(\gamma_0-\frac{2}{d+4}\right)^2}c_1 
& \mbox{ , } & c_1=\left(\sum_{j=1}^d\mu_j^2f^{(2)}_{jj}(x)\right)^2,\\
C_2\left(\xi\right)=\frac{\gamma_0}{2\left(\gamma_0-\frac{2}{d+4}\right)}c_2  
& \mbox{ , } & c_2=f\left(x\right)\int_{\mathbb{R}^d} K^2\left(z\right)dz,
\end{array}
\end{eqnarray*}
and the MSE of $f_n$ can be rewritten as
\begin{eqnarray*}
MSE=n^{-\frac{4}{d+4}}\delta_n^{\frac{4}{d+4}}\frac{d+4}{d^{\frac{d}{d+4}}
4^{\frac{d+6}{d+4}}}
\frac{\gamma_0^2}{\left(\gamma_0-\frac{2}{d+4}\right)^{\frac{2d+4}{d+4}}}
c_1^{\frac{d}{d+4}}c_2^{\frac{4}{d+4}}\left[1+o\left(1\right)\right].
\end{eqnarray*}
The function $x \mapsto x^2/\left(x-2/(d+4)\right)^{(2d+4)/(d+4)}$ reaching 
its minimum at the point $x=1$, to minimize the MSE of $f_n$, $\gamma_0$ 
must be chosen equal to $1$. Corollary~\ref{C:1} follows.

\paragraph{Acknowledgments} We are grateful to two Referees and an Associate Editor 
for their helpful comments, which have led to this substantially improved version of 
the paper.


\begin{thebibliography}{99}

\bibitem{arcones} Arcones, M.A. (1997),
The law of the iterated logarithm for a triangular array of empirical
processes,
{\em Electronic Journal of Probab.} {\bf 2}, 1-39


\bibitem{blu54}  Blum, J.R. (1954),
Multidimensional stochastic approximation methods,
{\em Ann. Math. Statist.,} {\bf 25},  737-744

\bibitem{Bo73}
 Bojanic, R. and Seneta, E. (1973), A unified theory of regularly varying sequences, \textit{Math}. Z., {\bf 134}, 91-106.

\bibitem{chen88}  Chen, H. (1988), 
Lower rate of convergence for locating a maximum of a function,
{\em Ann. Statist.,} {\bf 16},  1330-1334

\bibitem{2chen1999}  Chen, H.F., Duncan, T.E., and Pasik-Duncan, B. (1999), A Kiefer-Wolfowitz algorithm with randomized differences, {\em IEEE Trans. Automat. Control,} {\bf 44},  442-453


\bibitem{Dav73} Davies, H.L. (1973), Strong consistency of a sequential estimator of a probability density function. \textit{Bull. Math. Statist.}  \textbf{15},  49-54. 


\bibitem{deheu73} Deheuvels, P. (1973), \textit{Sur l'estimation s\'equentielle de la densit\'e, C. R. Acad. Sci. Paris Ser. A-B} \textbf{276}, 1119-1121.


\bibitem{De79} Devroye, L. (1979), On the pointwise and integral convergence of recursive kernel estimates of probability densities. \textit{Utilitas Math.}  \textbf{15},  113-128.

\bibitem{Di97}  Dippon, J. and Renz J. (1997),
Weighted means in stochastic approximation of minima,
{\em SIAM J. Control Optim.,} {\bf 35},  1811-1827

\bibitem{Di2003}  Dippon, J. (2003),
Accelerated randomized stochastic optimization,
{\em Ann. Statist.,} {\bf 31},  1260-1281

\bibitem{Du96} 
Duflo, M. (1996), \textit{Algorithmes stochastiques}, Collection Applications of mathematics, Springer.

\bibitem{Du97} 
Duflo, M. (1997), \textit{Random Iterative Models}, Collection Applications of mathematics, Springer.

\bibitem{Fabi67} Fabian, V. (1967),
Stochastic approximation of minima with improved asymptotic speed,
{\em Ann. Math. Statist.,} {\bf 38},  191-200

\bibitem{Ga73}
 Galambos, J. and Seneta, E. (1973), Regularly varying sequences, \textit{Proc. Amer. Math. Soc}., {\bf 41}, 110-116.

\bibitem{Hall80}  Hall, P. and Heyde, C.C. (1980), 
{\em Martingale limit theory and its application,} Academic Press, Inc., 
New York-London

\bibitem{hall} Hall, P. (1981), 
Laws of the iterated logarithm for nonparametric density estimators, 
{\em Z. Warsch. Verw. Gebiete} {\bf 56}, 47-61


\bibitem{hall1992} Hall, P. (1992), Effect of bias estimation on coverage accuracy of bootstrap confidence intervals for a probability density, {\em Ann. Statist.} {\bf 20}, 675-694


\bibitem{hall1994} Hall, P. and Patil, P. (1994), On the efficiency of on-line density estimators, 
{\em IEEE Trans. Inform. Theory} {\bf 40}, 1504-1512

\bibitem{Kie52}  Kiefer, J. and Wolfowitz, J. (1952),
Stochastic approximation of the maximum of a regression functions,
{\em Ann. Math. Statist.,} {\bf 23},  462-466

\bibitem{kushner1978} K{ushner}  H.J. and C{lark},  D.S. (1978) {\em Stochastic approximation methods for constrained and unconstrained systems,} Springer, New York

\bibitem{Men84} Menon, V.V., Prasad, B. and Singh, R.S. (1984), Non-parametric recursive estimates of a probability density function and its derivatives. \textit{Journal of Statistical Planning and inference}  \textbf{9},  73-82.

\bibitem{Mokka2007} Mokkadem, A. and  Pelletier, M. (2007a), 
{A companion for the Kiefer-Wolfowitz-Blum stochastic 
approximation algorithm,} \emph{Annals of Statistics}  {\bf 35}    no. 4, 1749--1772.


\bibitem{Mokka2006} Mokkadem, A. and  Pelletier, M. (2007b),
{Compact law of the iterated logarithm for matrix-normalized sums of random vectors,} 
{\em Teor. Veroyatn. Primen.}  {\bf 52}   




\bibitem{Pa62}
 Parzen, E. (1962). On Estimation of a Probability Density and Mode, \textit{Ann. Math. Statist}., {\bf 33}, 1065-1076.

\bibitem{Po90}  Polyak, B.T. and Tsybakov, A.B. (1990),
Optimal orders of accuracy for search algorithms of stochastic 
optimization,
{\em Problems Inform. Transmission,} {\bf 26},  126-133.

\bibitem{2RP73} R\'ev\'esz, P. (1973), Robbins-Monro procedure in a Hilbert space and its application in the theory of learning processes I. \textit{Studia Sci. Math. Hung.}, $\mathbf{8}$, 391-398.

\bibitem{RP77}
R\'ev\'esz, P. (1977), How to apply the method of stochastic approximation in the non-parametric estimation of a regression function, \textit{Math. Operationsforsch. Statist., Ser. Statistics}, {\bf 8}, 119-126.

\bibitem{Ro56}
 Rosenblatt, M. (1956), Remarks on Some Nonparametric Estimates of a Density Function, \textit{Ann. Math. Statist.}, {\bf 27}, 832-837.

\bibitem{Rou92} Roussas, G. (1992), Exact rates of almost sure convergence of a recursive kernel estimate of a probability density function: Application to regression and hazard rate estimate. \textit{J. of Nonparam. Statist.}  \textbf{3}, 171-195.

\bibitem{Ru82}  R{uppert}, D. (1982),
Almost sure approximations to the Robbins-Monro and Kiefer-Wolfowitz 
processes with dependent noise,
{\em Ann. of Probab.,} {\bf 10},  178-187.

\bibitem{spall1988} Spall, J.C. (1988), A stochastic approximation algorithm for large-dimensional systems in the Kiefer-Wolfowitz setting. In {\em Proc. Conference on Decision and Control}, 1544-1548. IEEE, New York.


\bibitem{Sp97} Spall, J.C. (1997),
A one-measurement form of simultaneous perturbation stochastic 
approximation, {\em  Automatica J. IFAC,} {\bf 33},  109-112.

\bibitem{Ts80} Tsybakov, A.B.  (1990),
Recurrent estimation of the mode of a multidimensional distribution,
{\em Problems Inform. Transmission} {\bf 26}, 31-37.
\bibitem{We79} Wegman, E.J. and Davies, H. I. (1979), Remarks on some recursive estimators of a probability density. \textit{Ann. Statist.}  \textbf{7},  316-327. 

\bibitem{Wol69} Wolverton, C.T. and Wagner, T.J. (1969), Asymptotically optimal discriminant functions for pattern classification. \textit{IEEE Trans. Inform. Theory}  \textbf{15},  258-265. 
\bibitem{Yama71} Yamato, H. (1971), Sequential estimation of a continuous probability density function and mode. \textit{Bull. Math. Satist.}  \textbf{14},  1-12. 
\end{thebibliography}
\end{document}